\newif\iflatexe
\latexefalse

\newif\ifdebug
\debugfalse

\newif\ifpreliminary
\preliminaryfalse

\iflatexe
\documentclass{article}
  \ifdebug
    \usepackage{latexsym,showlabels,mathptm,a4}
  \else
    \usepackage{latexsym,mathptm,a4}
  \fi
\else
  \ifdebug
    \documentstyle[a4,showlabels]{article}
  \else
    \documentstyle[a4]{article}
  \fi
\fi
\makeatletter
\@addtoreset{equation}{subsection}

\makeatother

\newcommand{\separate}{\medskip\noindent}

\def\newsection{ \separate
   \refstepcounter{subsection} 
   {\large\bf \thesubsection\kern.3em}
}
\def\mytheorem#1{
   \separate{\large\bf Theorem#1:\kern.3em}} 

\def\mylemma#1{
   \separate{\large\bf Lemma#1:\kern.3em}} 

\def\mycorollary#1{
   \separate{\large\bf Corollary#1:\kern.3em}} 

\def\myprop#1{
   \separate{\large\bf Proposition#1:\kern.3em}} 

\def\myremark#1{
   \separate{\large\bf Remark#1:\kern.3em}} 

\def\mydefinition#1{
   \separate{\large\bf Definition#1:\kern.3em}} 

\def\Proof{\separate\underline{Proof:}\kern1em}

\newcommand{\GANZ}{{\sf Z\hspace*{-0.4em}Z}}
\newcommand{\REELL}{{\setlength{\unitlength}{1em}
                     \begin{picture}(0.75,1)
                     \put(0,0){\line(0,1){0.69}}
                     \put(0,0){\sf R}
                     \end{picture}
                   }}
\newcommand{\KOMPLEX}{{\setlength{\unitlength}{1em}
                     \begin{picture}(0.7,1)
                     \put(0.34,0){\line(0,1){0.65}}
                     \put(0,0){\sf C}
                     \end{picture}
                   }}
\newcommand{\NATUR}{{\setlength{\unitlength}{1em}
                     \begin{picture}(0.75,1)
                     \put(0,0){\line(0,1){0.69}}
                     \put(0,0){\sf N}
                     \end{picture}
                   }}
\newcommand{\FGANZ}{\mbox{\tiny{\rm Z\hspace*{-0.45em}Z}}}
\newcommand{\FNATUR}{\mbox{\tiny{\setlength{\unitlength}{1em}
                     \begin{picture}(0.6,0.5)
                     \put(0,0){\line(0,1){0.48}}
                     \put(0,0){\rm N}
                     \end{picture}
                   }}}
\newcommand{\FKOMPLEX}{\mbox{\tiny{\setlength{\unitlength}{1em}
                               \begin{picture}(0.6,0.5)
                               \put(0.34,0){\line(0,1){0.47}}
                               \put(0,0){\rm C}
                               \end{picture}
                              }}}
\newcommand{\FREELL}{\mbox{\tiny{\setlength{\unitlength}{1em}
                     \begin{picture}(0.6,0.5)
                     \put(0,0){\line(0,1){0.48}}
                     \put(0,0){\rm R}
                     \end{picture}
                   }}}
\newcommand{\Brr}{{\mathchoice{\REELL}{\REELL}{\!\FREELL}{\!\FREELL}}}
\newcommand{\Bcc}{{\mathchoice{\KOMPLEX}{\KOMPLEX}{\!\FKOMPLEX}{
\!\FKOMPLEX}}}
\newcommand{\Bnn}{{\mathchoice{\NATUR}{\NATUR}{\FNATUR}{\FNATUR}}}
\newcommand{\Bii}{{\mathchoice{\GANZ}{\GANZ}{\FGANZ}{\FGANZ}}}
\newcommand{\CPE}{{\mathchoice{\Bcc{\rm P}_1}{\Bcc{\rm P}_1}{\Bcc\!\!
\mbox{\rm\tiny P}_1}{\Bcc\!\!\mbox{\rm\tiny P}_1}}}

\def\QED{\hfill$\Box$}
\def\inv{^{-1}}
\def\BEA{\begin{eqnarray}}
\def\EEA{\end{eqnarray}}
\def\BEQ{\begin{equation}}
\def\EEQ{\end{equation}}
\def\bref#1{(\ref{#1})}
\def\tmatrix#1#2#3#4{
    \left(\begin{array}{cc} #1 & #2 \\ #3 & #4 \end{array}\right)}
\def\LieSL{{\mbox{\bf SL}}}
\def\LieSO{{\mbox{\bf SO}}}

\def\LieSU{{\mbox{\bf SU}}}
\def\LieU{{\mbox{\bf U}}}
\def\LieB{{\mbox{\bf B}}}
\def\LieG{{\mbox{\bf G}}}
\def\Liesl{{\mbox{\bf sl}}}
\def\Liesu{{\mbox{\bf su}}}

\def\Lieg{{\mbox{\bf g}}}
\def\Liek{{\mbox{\bf k}}}
\def\Liep{{\mbox{\bf p}}}
\def\dprime{{\prime\prime}}

\def\vecr{{\mbox{\bf r}}}
\def\vecsigma{{\sigma\kern-2.25mm\sigma}}
\def\bpi{\bar{\pi}}
\def\bA{\bar{A}}

\def\tg{\tilde{g}}

\def\tv{\tilde{v}}
\def\tw{\tilde{w}}
\def\txi{\tilde{\xi}}

\def\tE{{\tilde{E}}}
\def\tF{{\tilde{F}}}
\def\tQ{{\tilde{Q}}}
\def\tR{{\tilde{R}}}
\def\tT{{\tilde{T}}}

\def\talpha{\tilde{\alpha}}
\def\tchi{\tilde{\chi}}
\def\tvarphi{\tilde{\varphi}}
\def\tepsilon{\tilde{\epsilon}}

\def\calpha{\check{\alpha}}
\def\cF{\check{F}}
\def\FA{{\cal A}}
\def\FC{{\cal C}}
\def\FD{{\cal D}}

\def\FI{{\cal I}}
\def\FL{{\cal L}}

\def\FP{{\cal P}}
\def\FR{{\cal R}}
\def\FRL{{{\cal R}\kern1mm_\lambda}}
\def\FS{{\cal S}}
\def\FT{{\cal T}}

\def\hpqt{{\overline{h_+}^\top}}

\def\Hqt{{\overline{H}^\top}}

\def\diff{{\rm d}}
\def\diffs{\diff s}
\def\diffz{\diff z}

\def\diffzbar{\diff\zbar}
\def\diffg{\diff g}
\def\diffp{\diff p}

\def\Ad{{\mbox{\rm Ad}}}

\def\diag{{\mbox{\rm diag}}}
\def\Ker{{\mbox{\rm Ker}}}

\def\Aut{{\mbox{\rm Aut}}}
\def\Iso{{\mbox{\rm Iso}}}
\def\OAff{{\mbox{\rm OAff}}}
\def\twospace{{\Brr^{\!\lower2pt\hbox{\mbox{\rm\scriptsize{2}}}}}}
\def\threespace{{\Brr^{\!\lower2pt\hbox{\mbox{\rm\scriptsize{3}}}}}}
\def\cstar{{\Bcc^{\!\lower2pt\hbox{\mbox{\scriptsize{$\ast$}}}}}}
\def\Imag{{\mbox{\rm Im}}}
\def\id{{\mbox{\rm id}}}

\def\const{{\mbox{\rm const}}}

\def\hxi{\hat{\xi}}

\def\hg{\hat{g}}
\def\hf{\hat{f}}
\def\hg{\hat{g}}
\def\hh{\hat{h}}
\def\hk{\hat{k}}
\def\hu{\hat{u}}
\def\hw{\hat{w}}
\def\hz{\hat{z}}

\def\hD{\hat{D}}
\def\hE{\hat{E}}
\def\hF{\hat{F}}

\def\hN{\hat{N}}

\def\htheta{\hat{\theta}}

\def\hU{\hat{U}}

\def\hPhi{\hat{\Phi}}
\def\hPsi{\hat{\Psi}}
\def\tg{\tilde{g}}

\def\tz{\tilde{z}}
\def\tQ{\tilde{Q}}
\def\tT{\tilde{T}}

\def\zbar{{\bar{z}}}

\def\Equer{{\overline{E}}}
\def\Squer{{\overline{S}}}

\def\qquer{{\overline{q}}}
\def\zquer{{\overline{z}}}
\def\alphaquer{{\overline{\alpha}}}
\def\betaquer{{\overline{\beta}}}

\def\gquer{{\overline{g}_-^\top}}
\def\hquer{{\overline{h}_+^\top}}

\def\rquer{{\overline{r}_+^\top}}
\def\hgquer{{\overline{\hg}_-^\top}}

\def\pder#1#2{{\partial #1\over\partial #2}}

\def\subcong{{\kern.3em\lower.4ex\hbox{$\stackrel{\displaystyle\subset}{\sim}$}\kern.3em}}
\def\supcong{{\kern.3em\lower.4ex\hbox{$\stackrel{\displaystyle\supset}{\sim}$}\kern.3em}}

\begin{document}

\renewcommand{\thefootnote}{\fnsymbol{footnote}}
  \vskip-6.94444pt
  \hbox to \textwidth{dg-ga/9603007\hfill March 19, 1996}
  \vskip6.94444pt
\begin{center}
{\LARGE On symmetries of constant mean curvature surfaces}

\vskip1cm
\begin{minipage}{6cm}
\begin{center}
J.~Dorfmeister\footnotemark[1]\\
Department of Mathematics\\ University of Kansas\\ Lawrence, KS 66045
\end{center}
\end{minipage}
\begin{minipage}{6cm}
\begin{center}
G.~Haak\footnotemark[2]\\ Fachbereich Mathematik\\ TU Berlin\\ D-10623 Berlin
\end{center}
\end{minipage}
\vspace{0.5cm}

\end{center}

\footnotetext[1]{partially supported by NSF Grant 
DMS-9205293 and Deutsche For\-schungs\-ge\-mein\-schaft}
\footnotetext[2]{supported by KITCS grant OSR-9255223 and
Sonderforschungsbereich 288}

\section{Introduction}\refstepcounter{subsection}\label{introduction}
\message{[introduction]}

The goal of this note is to start the investigation of conformal CMC-immersions
$\Psi:\FD\rightarrow\threespace$, $\FD$ an open, simply connected subset of
$\Bcc$, which allow for groups of spatial symmetries
$$
\Aut\Psi(\FD)=\{\tT\,\mbox{\rm proper Euclidean motion of $\threespace$}|
\tT\Psi(\FD)=\Psi(\FD)\}.
$$
More precisely (see the definition in Section~\ref{CMCautom}),
we consider a Riemann surface $M$ with universal covering
$\pi:\FD\rightarrow M$, and a conformal CMC-immersion
$\Phi:M\rightarrow\threespace$ with nonzero mean curvature,
such that $\Phi\circ\pi=\Psi$. Then we consider the groups
$$
\Aut\FD=\{g:\FD\rightarrow\FD\,\mbox{biholomorphic}\},
$$
$$
\Aut M=\{g:M\rightarrow M\,\mbox{biholomorphic}\},
$$
$$
\Aut_\pi\FD =\{g\in\Aut\FD| \;\mbox{\rm there exists}\; \hg\in\Aut M: \pi\circ g=\hg\circ\pi\},
$$
$$
\Aut_\Phi M =\{\hg\in\Aut M| \;\mbox{\rm there exists}\; \tT\in\Aut\Psi(\FD):
\Phi\circ\hg=\tT\circ\Phi\},
$$
and
$$
\Aut_\Psi\FD =\{g\in\Aut\FD| \;\mbox{\rm there exists}\; \tT\in\Aut\Psi(\FD):
\Psi\circ g=\tT\circ\Psi\}.
$$ 
There are many well known examples of CMC-surfaces with large spatial
symmetry groups. The
classic Delaunay surfaces (see~\cite{Ee:1}) 
have a nondiscrete group $\Aut\Psi(\FD)$
containing the group of all rotations around their generating axis.
Other examples are the Smyth surfaces \cite{Sm:2}, which were
visualized by D.~Lerner, I.~Sterling, C.~Gunn and U.~Pinkall. These
surfaces have an $m+2$-fold rotational symmetry in $\threespace$,
where the axis of rotation passes through the single umbilic of order
$m$. More recent is the large class of examples provided by 
K.~Gro\ss e-Brauckmann and K.~Polthier (see e.g.~\cite{GB:1,GBPo:1}) of singly,
doubly and triply periodic CMC-surfaces.

Other interesting classes of surfaces are the ones with a large group
$\Aut_\Psi\FD$. Examples for this are the compact
CMC-surfaces, whose Fuchsian or elementary group is contained in
$\Aut_\Psi\FD$. 

Yet another class of surfaces $(M,\Phi)$, 
for which $\Aut_\Psi\FD$ is interesting
are the surfaces with branchpoints. If one deletes the set $B\subset
M$ of those points in $M$
which are mapped by $\Phi$ to the branchpoints,
then such a surface can be constructed as an immersion $\hPhi$ of the
non-simply connected Riemann surface $M\setminus B$ into
$\threespace$. To get an immersion $\Psi$
of a simply connected domain $\FD$ into $\threespace$, which covers
$(M\setminus B,\hPhi)$, one can also apply the discussion of this paper as
will be shown in Section~\ref{branched}.

It is our goal to describe properties of the group $\Aut\Psi(\FD)$ in
terms of biholomorphic automorphisms of the Riemann surface $M$ or the simply
connected cover $\FD$, i.e., in terms of $\Aut_\Phi M$ or
$\Aut_\Psi\FD$.  To this end we investigate the relation
between these groups. This is done in Chapter~\ref{CMCautom}. After
defining in Section~\ref{setup},
what we mean by a CMC-immersions $(M,\Phi)$, we start
in Section~\ref{Riemannian} by listing some well known properties of
the groups $\Aut M$, $\Aut\FD$, $\Aut_\pi\FD$. These follow entirely
from the underlying Riemannian structure of $M$ and $\FD$. After
fixing the conventions for conformal CMC-immersions in
Section~\ref{DGL} we derive in Section~\ref{secfund} the
transformation properties of the metric and the Hopf differential
under an automorphism in $\Aut_\Psi\FD$. These will lead in
Section~\ref{Cautom} to some general restrictions on $\Aut_\Psi\FD$ in
the case $\FD=\Bcc$. In Section~\ref{uniqueness}--\ref{immersedautom}
we will introduce group homomorphisms $\bpi:\Aut_\pi\FD\rightarrow\Aut
M$, $\phi:\Aut_\Phi M\rightarrow\Aut\Psi(\FD)$ and
$\psi:\Aut_\Psi\FD\rightarrow\Aut\Psi(\FD)$. We will also prove, that,
in case $M$ with the metric induced by $\Psi$ is complete,
$\psi$ is surjective (Corollary~\ref{immersedautom}).  In
Section~\ref{groups} it will be shown that we furthermore can restrict
our investigation to those CMC-immersions
$\Phi:M\rightarrow\threespace$, for which $\phi$ is an isomorphism of
Lie groups.  In Sections~\ref{nondiscrete} and~\ref{regularsurfaces}
we will investigate, for which complete CMC-surfaces the group
$\Aut_\Psi\FD$ is nondiscrete. To this end we give a simple condition
on the image $\Psi(\FD)$ of $\Psi$ under which the group
$\Aut\Psi(\FD)$ is a closed Lie subgroup of $\OAff(\threespace)$. Here
we denote by $\OAff(\threespace)$ the group of proper
(i.e.,~orientation preserving) Euclidean motions of $\threespace$.
Using these results and the work of Smyth~\cite{Sm:2}, we will prove
in Section~\ref{Symsection}, that if
$\Aut\Psi(\FD)$ is closed, $\Aut\Psi(\FD)$ is nondiscrete, iff the surface
$\Psi(\FD)$ is isometric to a CMC-surface of revolution, i.e., a
Delaunay surface.  In Section~\ref{DelundSmyth} we illustrate the
discussion in Chapter~\ref{CMCautom}, using the examples of Delaunay
and Smyth. The groups $\Aut_\Phi M$, $\Aut_\Psi\FD$ and
$\Aut\Psi(\FD)$ are explicitly given for these examples.

In Chapter~\ref{DPWautom} we will recall the only known constructive
approach to general CMC-immersions, the DPW construction
(see~\cite{DoPeWu:1}). This will lead to the notion of the extended
frame for a conformal CMC-immersion (Section~\ref{Chistetig}). 
In Section~\ref{33}--\ref{loopinterpretation} we will derive
(Theorem~\ref{loopinterpretation}), how the elements of $\Aut_\Psi\FD$
act on the extended frame of a conformal CMC-immersion $\Psi$. This
will be used in Section~\ref{Symsection} to show,
how the elements of $\Aut_\Psi\FD$ act on the
whole associated family of a CMC-surface.
In Section~\ref{meropotentials} we define the meromorphic potential
$\xi$ of a CMC surface and
derive its transformation properties under $\Aut_\Psi\FD$.
This discussion will be used in Section~\ref{automorphism} to
characterize CMC-surfaces with symmetries 
in terms of their meromorphic
potentials (Proposition~\ref{automorphism}).

The discussion up to Section~\ref{automorphism} will be used in
Section~\ref{cyldressed} to prove (Theorem~\ref{cyldressed}),
that there are no nontrivial, translationally symmetric surfaces in
the dressing orbit of the cylinder. In the forthcoming
paper~\cite{DoHa:3}, the authors will use the results of these sections 
to give an alternative derivation of the classification of CMC-tori
in terms of algebro-geometric data~\cite{PiSt:1,Bo:1,ErKnTr:1,Ja:1}.

Chapter~\ref{DPWautom} closes with an investigation of the
coefficients of the meromorphic potential and their transformation
properties. We will give (Theorem~\ref{complete})
a necessary condition on the meromorphic potential for the surface
$(M,\Phi)$ to be complete. 

In Chapter~\ref{unitary} we start a closer investigation of the
transformation properties we derived for the meromorphic potential. In
Section~\ref{summaryth} we
will prove necessary and sufficient conditions on the coefficient 
functions of the meromorphic potential in order to give
surfaces with symmetries under the DPW construction. 
Unfortunately, these conditions 
still involve the Birkhoff splitting of the extended
frames. Ne\-ver\-the\-less, we expect them to be suitable starting points for
the investigation and construction of many concrete examples of
surfaces with symmetries and, in particular, of compact CMC-surfaces.

We finish the paper in Chapter~\ref{examples} with a closer look at
the meromorphic potential of Smyth surfaces and the treatment of a
CMC-immersion with a branchpoint.

The results of sections~\ref{33}, \ref{loopinterpretation} and
\ref{meropotentials} are based on discussions of one of the authors
(J.D.) with Franz Pedit.

\section{Automorphisms of CMC-surfaces} \label{CMCautom}\message{[CMCautom]}

\newsection \label{setup}
Before we can start the investigation of CMC-immersions $(M,\Phi)$ we
have to define, what we mean by a CMC-immersion, if $M$ is not a
domain in $\twospace$. We will restrict our investigations to the case
of nonzero mean curvature, i.e.\ we will exclude the special case of
minimal surfaces.

\mydefinition{}
Let $M$ be a connected $C^2$-manifold and let $\Phi:M\rightarrow\threespace$ be
an immersion of type $C^2$. $\Phi$ is called a CMC-immersion, if
there exists an atlas of $M$, s.t.\ every chart $(U,\varphi)$ 
in this atlas defines a $C^2$-surface
$\Phi\circ\varphi\inv:\varphi(U)\rightarrow\threespace$ 
with nonzero constant mean curvature.

\separate
That this definition makes sense is shown by the following

\mytheorem{} {\em
Let $M$ be a connected $C^2$-manifold for which there exists a $C^2$-immersion
$\Phi:M\rightarrow\threespace$, s.t.\ $(M,\Phi)$ is a
CMC-immersion. Then there exists an atlas of $M$, s.t.\ the
mean curvature is globally constant. In particular, $M$ can be oriented and
the mean curvature depends only on the chosen orientation of $M$.
}

\Proof
Let $\FA=\{(U_\alpha,\varphi_\alpha),\alpha\in\FI\}$, $\FI$ some index
set, be an atlas of $M$, s.t.\ for every chart
$\Phi\circ\varphi_\alpha\inv:\varphi_\alpha(U_\alpha)\rightarrow\threespace$
defines a $C^2$-surface of constant mean curvature $H_\alpha\neq0$.
By changing, if necessary, the
orientation of the surfaces $\Phi\circ\varphi_\alpha\inv$, we can
assume that $H_\alpha>0$ for all $\alpha\in\FI$.
Let $U_\alpha$ and $U_\beta$, $\alpha\neq\beta$, be s.t.\ $U_\alpha\cap
U_\beta\neq\{\}$. Then, on $U_\alpha\cap U_\beta$,
$\Phi\circ\varphi_\alpha\inv$ and $\Phi\circ\varphi_\beta\inv$ define
two surfaces which differ only by a $C^2$-reparametrization
$\varphi_\beta\circ\varphi_\alpha\inv$.
The absolute value of the mean curvature is
independent of reparametrizations. Since $H_\alpha$ and $H_\beta$
are positive, we get $H_\alpha=H_\beta$ and
$\varphi_\beta\circ\varphi_\alpha\inv$ has positive Jacobian.
{}From this it follows that $M$ is orientable. Since $M$ is also connected, 
we get that all $H_\alpha$ coincide. Therefore, the
mean curvature defined by $\FA$ on $M$ is globally constant.

By the arguments above, for a chosen orientation on $M$, 
the absolute value of the mean curvature $H$ is
constant on $M$ and does not depend on the chosen atlas, as long as
this atlas defines the same orientation. Changing the
orientation of the atlas changes the sign of $H$, which finishes the proof.
\QED

\newsection
We will not work with the fairly general Definition~\ref{setup}.
Instead we use Lichtenstein's theorem to turn $M$ into
a Riemann surface (see also~\cite[Theorem~5.13]{Spr:1})
and $\Phi$ into a conformal immersion.

\mytheorem{} {\em Let $(M,\Phi)$ be a CMC-immersion. Then there exists a
conformal structure on $M$, s.t.\ $M$ becomes a Riemann surface
and $\Phi$ becomes a conformal
CMC-immersion.
}

\Proof
Let $\FA=\{(U_\alpha,\varphi_\alpha),\alpha\in\FI\}$ be an oriented atlas of
$M$. W.l.o.g.\ we can choose $\FA$ s.t. $\Phi|_{U_\alpha}$
is injective for all $\alpha\in\FI$. Then, by Lichtenstein's theorem (see
e.g.~\cite[Section~I.1.4]{DiHiKuWo:1}), for every $\alpha\in\FI$,
there exists an orientation preserving diffeomorphism
$q_\alpha:V_\alpha\rightarrow\varphi(U_\alpha)$, $V_\alpha$ a
domain in $\twospace$, s.t.\ $\Phi\circ\varphi\inv_\alpha\circ
q_\alpha$ is a conformal parametrization of the surface defined by
$\Phi\circ\varphi\inv_\alpha$. Then
$\{(U_\alpha,\tvarphi_\alpha=q_\alpha\inv\circ\varphi_\alpha),\alpha\in\FI\}$
defines an atlas on $M$, for which every transfer function is holomorphic,
which finishes the proof.
\QED

\separate
In this paper we will therefore restrict ourselves to conformal
CMC-immersions $\Phi:M\rightarrow\threespace$, where $M$ is a
Riemann surface. I.e., if we write ``{\em $(M,\Phi)$ is a 
CMC-immersion}'', we always mean ``{\em $M$ is a Riemann surface and
$\Phi$ is a conformal CMC-immersion}''.

\newsection \label{Riemannian}
In this section we strip the CMC-surface $M$ of its metric and leave only
the complex structure. We will recollect some well known facts
about Riemann surfaces.

Up to conformal equivalence the only Riemann surfaces, 
which are simply connected are the
sphere $\CPE\cong\Bcc\cup\{\infty\}$, the complex plane $\Bcc$ and the upper
half plane $\Delta=\{z\in\Bcc | \Imag(z)>0\}$, which is conformally
equivalent to the open unit disk. Each of these surfaces is equipped
with its standard complex structure.

For convenience of language we will not distinguish between Fuchsian
groups (as in the case of $\FD=\Delta$) and elementary groups (as in
the case of $\FD=\Bcc$). For a more thorough treatment of the
uniformization of Riemann surfaces see e.g.\ \cite{FaKr:1}.

Every Riemann surface $M$ can be represented as the quotient of one of
these three Riemann surfaces by a freely acting Fuchsian group $\Gamma$ of
biholomorphic automorphisms, i.e., we
may write $M=\Gamma\setminus\FD$, where $\FD$ is the simply connected
cover of $M$.

If $\pi:\FD\rightarrow M$ is the covering map, then $\Gamma$ is also
the covering group of $\pi$. Then by~\cite[V.4.5]{FaKr:1}
\BEQ
\Aut M=N(\Gamma)/\Gamma,
\EEQ
where $N(\Gamma)$ is the normalizer of $\Gamma$ in $\Aut\FD$.
Since $g\in\Aut_\pi\FD$ is equivalent to 
\BEQ
\pi\circ g\circ\gamma=\pi\circ g
\EEQ
for all $\gamma\in\Gamma$, we have $g\circ\gamma=\gamma_1\circ g$ for
some $\gamma_1\in\Gamma$, and thus
\BEQ \label{AutpiFD}
\Aut_\pi\FD=N(\Gamma).
\EEQ
Therefore,
\BEQ \label{AutMcong}
\Aut M\cong\Aut_\pi\FD/\Gamma,
\EEQ
where $\Aut_\pi\FD$ is a closed subgroup 
of $\Aut\FD$, and $\Gamma$ is a normal subgroup of $\Aut_\pi\FD$.

The following is also well known (see e.g.\ \cite[IV.5,IV.6,V.4]{FaKr:1}):

\mylemma{} {\em
\begin{description}
\item[a)] The group $\Gamma$ is discrete, i.e., either finite or countable, and
consists of conformal (biholomorphic) automorphisms of $\FD$, which
act fixed point free on $\FD$.
\item[b)] If $\FD=\CPE$, then $\Gamma$ is trivial and $M$ is the
sphere.
\item[c)] If $\FD=\Bcc$, then $\Gamma$ is abelian and $M$ is either the plane,
the cylinder or a torus.
\item[d)] For a Riemann surface $M=\Gamma\setminus\FD$ the
following are equivalent:
\begin{itemize}
\item The Fuchsian group is abelian.
\item The Lie group $\Aut M$ of conformal automorphisms
of $M$ is nondiscrete.
\end{itemize}
Surfaces of this kind are called exceptional surfaces.
\end{description}
}

\newsection \label{DGL}
Now we also take into account the metric of the CMC-surface. 

Let $(M,\Phi)$ be a CMC-immersion. 
The Riemann surface $M$ is endowed with the induced
metric. As usual, $M$ is called (geodesically) complete as a manifold with
metric, if each geodesic can be extended
to a curve parametrized over $\Brr$. It is not possible to reach the
boundary of a complete manifold by going along a curve of finite length.

We are interested in orientation preserving isometries of the surface 
onto itself. They are automatically biholomorphic automorphisms, 
but may have fixed points.

We consider the commutative diagram
\BEQ \label{Riemanndiagram}
\setlength{\unitlength}{0.00083300in}%
\begin{picture}(1265,1215)(4316,-2269)
\thinlines
\put(4801,-1261){\vector(-1,-2){300}}
\put(4951,-1261){\vector( 1,-2){300}}
\put(4501,-2086){\makebox(0,0)[b]{$M$}}
\put(4651,-2011){\vector( 1, 0){450}}
\put(5251,-2086){\makebox(0,0)[b]{$\threespace$}}
\put(5251,-1561){\makebox(0,0)[b]{$\Psi$}}
\put(4801,-2236){\makebox(0,0)[b]{$\Phi$}}
\put(4876,-1186){\makebox(0,0)[b]{$\FD$}}
\put(4501,-1561){\makebox(0,0)[b]{$\pi$}}
\end{picture}
\EEQ
where $M$ is a Riemann surface and $\Phi$ is a CMC-immersion
of $M$.
Moreover, $\pi$ is the universal covering map of $M$ and
$\Psi=\Phi\circ\pi$. Remember, that we always assume,
that $\Phi$ and therefore also $\Psi$ is a conformal immersion.

\separate
We would like to remind the reader of the following
well known result (see e.g.\ the Appendix of \cite{DoHa:1}):

\mytheorem{} {\em
Let $\Psi:\FD\subset\Bcc\rightarrow\threespace$ be a
conformal immersion with metric
\BEQ \label{udef}
\diffs^2=\frac{1}{2}e^u\diffz\diffzbar
\EEQ
where $u=u(z,\zquer):\FD\rightarrow\Brr$. Define the function
$E:\FD\rightarrow\Bcc$ by
\BEQ \label{Hopfdef}
E=\langle\Psi_{zz},N\rangle
\EEQ
where $\langle\cdot,\cdot\rangle$ is the standard scalar product in
$\threespace$ and
\BEQ
N=\frac{\Psi_z\times\Psi_\zquer}{|\Psi_z\times\Psi_\zquer|}
\EEQ
is the Gau\ss\ map of $(\FD,\Psi)$.
Then the second fundamental form of $(\FD,\Psi)$ has in real
coordinates $x=\frac{1}{2}(z+\zquer)$, $y=\frac{1}{2i}(z-\zquer)$, the form
\BEQ \label{secfundform}
II=\frac{1}{2}\tmatrix{(E+\Equer)+He^u}{i(E-\Equer)}{i(E-\Equer)}{
-(E+\Equer)+He^u},
\EEQ
where $H=2e^{-u}\langle\Psi_{z\zquer},N\rangle$ is the mean curvature.
The Gau\ss-Codazzi equations take the form
\BEA
u_{z\zquer}+{1\over2}e^uH^2-2e^{-u}|E|^2 & = & 0, \\
E_\zquer={1\over2}e^uH_z.\label{Ezquer}
\EEA
The Gau\ss\ curvature is given in terms of $u$, $E$ and $H$ by
\BEQ \label{Gausscurvature}
K=H^2-4|E|^2e^{-2u}.
\EEQ
}

\separate From this the following corollary follows immediately

\mycorollary{} {\em
Let $\FD$ be the open unit disk or the complex plane.
\begin{description}
\item[1.] Let $\Psi:\FD\rightarrow\threespace$ be a CMC-immersion. Let 
$u$ and $E$ be defined as above.
Then $E\diffz^2$ is a holomorphic two-form on $\FD$ and the 
function $u$ satisfies the partial differential equation
\BEQ \label{gensinh}
u_{z\zbar}+\frac{H^2}{2}e^u-2|E|^2e^{-u}=0,
\EEQ
where $H$ is the mean curvature.
\item[2.] Let $E\diffz^2$ be a holomorphic differential on $\FD$ and
$u(z,\zquer)$ be a real valued function on $\FD$. Assume, that $u$ and $E$ satisfy
Eq.~\bref{gensinh} for some constant $H$. Then there exists a
conformal immersion $\Psi:\FD\rightarrow\threespace$ with constant
mean curvature $H$, s.t.\ $u$ and $E$ are
given by Eqs.~\bref{udef} and~\bref{Hopfdef},
respectively. The conformal immersion $\Psi$ is unique up to a
proper Euclidean transformation.
\end{description}
}

\Proof
1.\ follows from Eqs.~\bref{gensinh} and~\bref{Ezquer}.

2.\ For $u$ and $E$ the Gauss-Codazzi equations~\bref{gensinh} 
and~\bref{Ezquer}
are satisfied. The functions $u$ and $E$ determine, by~\bref{udef}
and~\bref{secfundform}, the first and second fundamental form
uniquely. The existence and the uniqueness of $\Psi$, up to proper
Euclidean motion,
follows therefore from the fundamental theorem of surface theory 
(see e.g.~\cite[III.2.3]{Sp:1}).
\QED

\newsection \label{associatedfamily}
For a CMC-immersion the holomorphic 
two-form $E\diffz^2$ is called the
Hopf differential. It is identically zero iff the surface
$(\FD,\Psi)$ is part of a sphere. Its zeroes are mapped to the
umbilical points of the surface.
The cylinder of mean curvature $H$ can be represented by
the ``vacuum solution'' $E\equiv\const\neq0$, 
$e^{2u}=\frac{4}{H^2}|E|^2\equiv\const$, of Eq.~\bref{gensinh}. 
As a special case of Corollary~\ref{DGL} we get the

\myprop{} {\em
If $\Psi:\FD\rightarrow\threespace$ is a complete CMC-immersion with
$u\equiv\const$ and $E\equiv\const\neq0$, 
then, up to orientation, $\FD=\Bcc$ and $\Psi(\FD)$ 
is a cylinder of mean curvature $H=2e^{-u}|E|$.
}

\Proof Since, by assumption, $\FD$ with the
constant metric $\diffs^2=\frac{1}{2}e^u\diffz\diffzbar$ is complete,
we have $\FD=\Bcc$. If we define $H=2e^{-u}|E|$ or $H=-2e^{-u}|E|$,
then $u$ and $E$ satisfy Eq.~\bref{gensinh}. Moreover, these are the
only two choices of $H$, for which Eq.~\bref{gensinh} is satisfied.
Therefore, $\Psi$ has the same metric and Hopf differential as a
cylinder of mean curvature $H=\pm2e^{-u}|E|$, where the sign depends
only on the choice of the orientation.
The rest of the proposition then follows from the
second part of Corollary~\ref{DGL}.
\QED

\separate
Since, by a scaling in $\threespace$ and the choice of an orientation
of $M$, $H$ can be fixed to any numerical 
value in $\Brr\setminus\{0\}$, we will from now on set
\BEQ
H=-\frac{1}{2}.
\EEQ
Let us fix a CMC-immersion $(\FD,\Psi)$ with Hopf 
differential $E$ and let $u$ be defined by Eq.~\bref{udef}.
Then $u$ solves Eq.~\bref{gensinh} for the given $E$.

It is clear, that $u$ still satisfies
Eq.~\bref{gensinh} if we replace $E$ by $\lambda E$, $|\lambda|=1$.
In this way, by Corollary~\ref{DGL}, we get for every
$\lambda\in S^1$ a CMC-immersion
$\Psi_\lambda:\FD\rightarrow\threespace$ 
with constant mean curvature $H$. The surface $(\FD,\Psi_\lambda)$ is
uniquely determined up to a proper Euclidean motion.
The family $\{\Psi_\lambda,\lambda\in S^1\}$ of CMC-immersions 
is called the {\em associated family} of $(\FD,\Psi=\Psi_{\lambda=1})$.
Note, that up to a proper Euclidean motion,
\BEQ \label{assgroup}
(\Psi_\lambda)_\mu=\Psi_{\lambda\mu},\kern1cm\mbox{\rm for
$\lambda,\mu\in S^1$.} 
\EEQ

Let $\Psi_1$ and $\Psi_2$ be two conformal immersions from $\FD\subset\Bcc$ to
$\threespace$ with metrics $\diff s_1^2=\frac{1}{2}e^{u_1}\diffz\diffzbar$
and $\diff s_2^2=\frac{1}{2}e^{u_2}\diffz\diffzbar$.
We call $\Psi_1$ and $\Psi_2$ {\em isometric} iff there is a conformal
automorphism $g:\FD\rightarrow\FD$, s.t.\ $\Psi_1\circ g$ and $\Psi_2$ have the
same metric, i.e., iff $e^{(u_1\circ g)}|g^\prime|^2=e^{u_2}$.

\mylemma{} {\em
Let $\Psi_1$ and $\Psi_2$ be two CMC-immersions from
$\FD\subset\Bcc$ to $\threespace$ with the same mean curvature.
Then $\Psi_1$ and $\Psi_2$ are isometric iff there is a biholomorphic
automorphism $g$, s.t.\ $\Psi_1\circ g$ and $\Psi_2$ are in the same
associated family.
}

\Proof 
If $g$ is a conformal automorphism of $\FD$, then the metric of
$\hPsi=\Psi\circ g$ is given by $\frac{1}{2}e^{\hu}\diffz\diffzbar$ with
$e^{\hu}=e^{(u_1\circ g)}|g^\prime|^2$, 
and the Hopf differential of $\Psi\circ g$ is $\hE=(E\circ g)(g^\prime)^2$.

Let $\Psi_1$ and $\Psi_2$ be two CMC-immersions with the same mean
curvature, and let $g$ be a
conformal automorphism, s.t.\ $\hPsi=\Psi_1\circ g$ and $\Psi_2$ have
the same metric, i.e.\ $e^{\hu}=e^{(u_1\circ g)}|g^\prime|^2=e^{u_2}$.
Since $\Psi_1$ and $\Psi_2$ have the same mean curvature, $g$ is
orientation preserving, i.e.\ $g$ is a biholomorphic automorphism of
$\FD$. In addition, 
since the Gau\ss\ curvature is invariant under isometries, it follows
from~Eq.~\bref{Gausscurvature} that
\BEQ 
|(E_1\circ g)(z)|\cdot|g^\prime(z)|^2=|E_2(z)|,
\EEQ
i.e.\ $\hE(z)=(E_1\circ g)(z)(g^\prime(z))^2=E_2(z)e^{i\theta(z)}$, where
$\theta(z):\FD\rightarrow\Brr$. Since $\hE$ and $E_2$ are holomorphic
we get that $\theta$ is constant. Consequently, 
$\hPsi$ and $\Psi_2$ are in the same associated family.

Conversely, assume that there exists a biholomorphic automorphism $g$,
s.t.\ $\Psi_1\circ g$ and $\Psi_2$ are in the same associated
family. Then the functions $e^{\hu}=e^{(u_1\circ g)}|g^\prime|$ and $e^{u_2}$
coincide and therefore $\Psi_1$ and $\Psi_2$ are isometric.
\QED

\separate It is clear from the definition of a conformal CMC-immersion, that
with $(M,\Phi)$, $M$ a Riemann surface, also $(M,\Phi\circ\hg)$,
$\hg\in\Aut M$, is a conformal CMC-immersion.
If $(\FD,\Psi)$ is a CMC-immersion from $\FD\subset\Bcc$ to $\threespace$,
then we call $g\in\Aut\FD$ a {\em self-isometry} of $(\FD,\Psi)$ if
the CMC-immersions $\Psi\circ g$ and $\Psi$ have the same
metric. If $(M,\Phi)$ is a CMC-immersion with universal cover
$(\FD,\Psi)$, $\pi:\FD\rightarrow M$, 
then $\hg\in\Aut M$ will be called a self-isometry of
$(M,\Phi)$ if there exists a self-isometry $g$ of $(\FD,\Psi)$, s.t.\
$\pi\circ g=\hg\circ\pi$.
We denote by $\Iso_\Psi\FD\subset\Aut\FD$ and $\Iso_\Phi M\subset\Aut
M$ the groups of self-isometries of $(\FD,\Psi)$ and $(M,\Phi)$, respectively. 
Note, that by definition, a self-isometry is orientation preserving.
{}From the proof of Lemma~\ref{associatedfamily} we get

\mycorollary{ 1} {\em Let $\Psi:D\rightarrow\threespace$ be a 
CMC-immersion with Hopf differential $E\diffz^2$ and define the
function $u$ as in Eq.~\bref{udef}. Let $g\in\Aut\FD$. Then the
following are equivalent:
\begin{description}
\item[1.] The automorphism $g$ is in $\Iso_\Psi\FD$.
\item[2.] The function $u$ transforms under $g$ as
\BEQ \label{uselfisomtrafo}
e^{(u\circ g)(z,\zquer)}|g^\prime(z)|^2=e^{u(z,\zquer)}.
\EEQ
\end{description}
If $g\in\Iso_\Psi\FD$, then $E$ transforms under $g$ as
\BEQ \label{Eselfisomtrafo}
|(E\circ g)(z)|\cdot|g^\prime(z)|^2=|E(z)|.
\EEQ
}

Since in Chapter~\ref{DPWautom} we will consider only complete CMC
immersions, we also want to state the following

\mycorollary{ 2} {\em 
Let $(\FD,\Psi)$ be a CMC-immersion. 
If $(\FD,\Psi)$ is complete (w.r.t.\ the induced metric),
then all elements $(\FD,\Psi_\lambda)$, $\lambda\in S^1$, of
its associated family are complete. If $(\FD,\Psi)$ is not complete,
then no element of its associated family is complete.
}

\Proof Since all surfaces in the associated
family share the same metric, they are either all complete or none of
them is complete.
\QED

\newsection \label{secfund}
We define $\OAff(\threespace)$ to be the group of proper Euclidean 
motions in $\threespace$.
It will be convenient at times to decompose an element of
$\tT\in\OAff(\threespace)$ into a rotational and a translational part:
\BEQ \label{congdecomp}
\tT v=R_{\tT}v+t_{\tT},\kern1cm v\in\threespace.
\EEQ
We will also write $\tT=(R_{\tT},t_{\tT})$.

\mydefinition{}
As already mentioned in the introduction we define
\BEQ \label{eqa}
\Aut_\pi\FD =\{g\in\Aut\FD| \;\mbox{\rm there exists}\; \hg\in\Aut M: \pi\circ g=\hg\circ\pi\},
\EEQ
\BEQ \label{eqb}
\Aut_\Phi M=\{\hg\in\Aut M|\;\mbox{\rm there exists}\; \tT\in\Aut\Psi(\FD):
\Phi\circ\hg=\tT\circ\Phi\},
\EEQ
\BEQ \label{eqc}
\Aut_\Psi\FD=\{g\in\Aut\FD|\;\mbox{\rm there exists}\; \tT\in\Aut\Psi(\FD):
\Psi\circ g=\tT\circ\Psi\},
\EEQ
and
\BEQ \label{eqd}
\Aut\Psi(\FD)=\{\tT\in\OAff(\threespace)|
\tT\Psi(\FD)=\Psi(\FD)\}.
\EEQ

\mylemma{} {\em
Let $\Psi:\FD\rightarrow\threespace$ be a CMC-immersion with Hopf
differential $E\diffz^2$ and define the
function $u$ as in Eq.~\bref{udef}. Let $g\in\Aut\FD$. 
Then the following are equivalent:
\begin{description}
\item[1.] The automorphism $g$ is in $\Aut_\Psi\FD$.
\item[2.] The functions $u$ and $E$ transform under $g$ as
\BEQ \label{utrafo}
e^{(u\circ g)(z,\zquer)}|g^\prime(z)|^2=e^{u(z,\zquer)},
\EEQ
\BEQ \label{Etrafo1}
(E\circ g)(z)(g^\prime(z))^2=E(z).
\EEQ
\end{description}
}

\Proof Let us define the immersion
$\Psi_1=\Psi\circ g$. Then $\Psi_1:\FD\rightarrow\threespace$
is also a CMC-immersion. By the definition of $u$ we have
\BEQ \label{u1eq}
e^{u_1(z,\zquer)}=e^{(u\circ g)(z,\zquer)}|g^\prime(z)|^2.
\EEQ
Since the Hopf differential is a holomorphic two-form we get
\BEQ \label{E1eq}
E_1(z)=(E\circ g)(z)(g^\prime(z))^2.
\EEQ
We have $g\in\Aut_\Psi\FD$ iff $\Psi_1$ and
$\Psi$ give the same surface up to a proper Euclidean motion.
By the fundamental theorem of surface theory this is the case iff
both surfaces have the same first and second fundamental form, which
by Eqs.~\bref{udef} and~\bref{secfundform} is equivalent to
$E_1=E$ and $u_1=u$. This, together with Eq.~\bref{u1eq} and
Eq.~\bref{E1eq}, proves the lemma.
\QED

\mycorollary{} {\em The elements of $\Aut_\Psi\FD$ act as
self-isometries of $(\FD,\Psi)$, i.e.\ $\Aut_\Psi\FD\subset\Iso_\Psi\FD$.
}

\Proof By Eq.~\bref{utrafo} and the definition
\bref{udef} of $u$, the metric $\diffs^2$
is invariant under $g\in\Aut_\Psi\FD$. \hbox{\kern3cm}
\QED

\newsection \label{Cautom}
We will draw some further conclusions from the automorphicity of the
Hopf differential, Eq.~\bref{Etrafo1}. We recall that CMC-surfaces with Hopf
differential identically zero are part of a round sphere. Such
surfaces will be called {\em spherical}.

\myprop{} {\em Let $(M,\Phi)$ be a CMC-surface with simply
connected cover $(\FD,\Psi)$, $\FD=\Bcc$. Then either $E\equiv0$ or 
the group $\Iso_\Psi\FD\subset\OAff(\threespace)$ 
of self-isometries of $(\FD,\Psi)$ consists only of
rigid motions of the plane, i.e., every $g\in\Iso_\Psi\FD$ can be
written as $g:z\mapsto az+b$, with $a,b\in\Bcc$ and $|a|=1$.
}

\Proof Let us assume, that there exists 
an automorphism $g$ in $\Iso_\Psi\FD\subset\Aut\Bcc$, which is of the form
\BEQ \label{grep}
g(z)=az+b
\EEQ 
with $a,b$ being complex constants and $|a|\neq1$.

Case I: If $b\neq0$ then we can, by a biholomorphic change of coordinates
\BEQ
z\mapsto\tz=\frac{a-1}{b}z+1,
\EEQ
turn $g$ into a scaling with rotation $g(\tz)=a\tz$. Let us also
define $\tE:\FD\rightarrow\Bcc$ by $\tE\diff\tz^2=E\diffz^2$, then
\BEQ
\tE(\tz)=\frac{b^2}{(a-1)^2}E(\frac{b}{a-1}(z-1)).
\EEQ

The Hopf differential $\tE$ transforms under $\tg$ according to 
\bref{Eselfisomtrafo}.
We therefore get for all $n\in\Bii$:
\BEQ \label{Edecrease}
|\tE(a^n\tz)|\cdot|a^{2n}|=|\tE(\tz)|.
\EEQ
For $|a|>1$ this implies that the absolute value of $\tE$ is decreasing
from the fixed point $\tz=0$ of $g$ in all
directions in the $\tz$-plane to zero. 
Since $\tE$ is holomorphic in $\tz$, this gives $\tE\equiv0$ and
therefore $E\equiv0$.

If $|a|<1$, consider the inverse $g\inv\in\Iso_\Psi\FD$:
\BEQ
g\inv(z)={1\over a}z-{b\over a}.
\EEQ
Since $\left|{1\over a}\right|>1$
we can use the first part of the proof again.

Case II: If $b=0$ then Eq.~\bref{Eselfisomtrafo} gives directly
\BEQ \label{Edecrease2}
|E(a^nz)|=|a^{-2n}||E(z)|.
\EEQ
We can therefore argue in the same way as in the first case.
\QED

\separate
The following Lemma will also be important:

\mylemma{} {\em
Let $(M,\Phi)$ be a complete, nonspherical
CMC-surface with conformal covering immersion $(\FD,\Psi)$ and
$\FD=\Bcc$. Then the following holds
\begin{description}
\item[1.] If $\Iso_\Psi\FD$ contains the one-parameter group of
rotations $\FR$ around a fixed point and also a translation, then $\Phi(M)$
is a cylinder.
\item[2.] If $\Iso_\Psi\FD$ contains a one-parameter group of
translations $\FT$ and also a rotation $R$, s.t.\ $R^2\neq\id$, then $\Phi(M)$
is a cylinder.
\end{description}
}

\Proof 1.\ Let $u$ be defined by
Eq.~\bref{udef}. W.l.o.g.\ we can 
choose the center of the rotation group $\FR$ as
$z=0$ and the translation as $g:z\mapsto z+1$.
Then $\FR=\{r_\phi|r_\phi(z)=e^{i\phi}z, \phi\in[0,2\pi)\}$.
By Eqs.~\bref{uselfisomtrafo} and~\bref{Eselfisomtrafo} we have
\BEQ
e^{u(z+1)}=e^{u(z)},\kern3cm |E(z+1)|=|E(z)|
\EEQ
\BEQ
|e^{u(e^{i\phi}z)}|=|e^{u(z)}|,\kern3cm |E(r_\phi(z))|=|E(e^{i\phi}z)|=|E(z)|
\EEQ
for all $z\in\Bcc$.
Therefore, $e^u$ and $|E|$ are constant
on the whole orbit of $z=0$ under the group
$\FR\times\FT$, where $\FT$ is the group generated by $z\mapsto z+1$.
The translation $z\mapsto z+1$ takes the origin to $z=1$, which is
mapped to the whole unit circle by the rotation group $\FR$.
The group $\FT$ takes the unit circle into
a connected set reaching from $0$ to $\infty$. This set therefore
contains for each
$r\in\Brr$ a complex number $z_r$ with $|z_r|=r$. Thus, the orbit
of this set under $\FR$ is the whole complex plane.
This shows, that $\FR\times\FT$ acts transitively on $\Bcc$ and 
$e^u$ and $|E|$ are constant functions. Since $u$ is real valued and $E$
is holomorphic, we get that $u$ and $E$ are constant functions.
The image $\Phi(M)$ is then, by Proposition~\ref{associatedfamily}, a cylinder.

2.\ W.l.o.g.\ we choose the group $\FT$ as the group of
translations along the real axes $\{z\mapsto z+r|r\in\Brr\}$ and $R$
as a rotation around the origin $R(z)=e^{i\phi}$ with $\phi\neq
m\pi$, $m\in\Bii$. Let us denote the group generated by $R$ as $\FR$.
The orbit of the origin $z=0$ under the group $\FR\times\FT$ contains
the straight line $\{z|z=re^{i\phi},r\in\Brr\}$, which does not 
coincide with the
real line. Therefore, the action of $\FT$ takes this line into the
whole complex plane and the group $\FR\times\FT$ acts transitively on $\Bcc$.
As in the proof of the first part of the lemma
this shows, that $E$ and $u$ are constant and that $\Phi(M)$ is a cylinder.
\QED

{}From the results of this section we can draw the following conclusion
for $\Iso_\Psi(\FD)$, if $\FD=\Bcc$:

\mytheorem{} {\em
Let $(M,\Phi)$ be a complete, nonspherical CMC-surface
with universal covering
immersion $(\FD,\Psi)$ and $\FD=\Bcc$.
Let $\Iso_\Psi\FD\subset\Aut\FD$ be the
group of self-isometries of $(\FD,\Psi)$.

1. If $\Iso_\Psi\FD$ contains the one-parameter group 
$\FR$ of rotations around a fixed point $z_0\in\Bcc$, then 
\begin{itemize}
\item either $\Iso_\Psi\FD=\FR$
\item or $\Phi(M)$ is a cylinder.
\end{itemize}

2. If $\Iso_\Psi\FD$ contains a one-parameter group $\FT$ of translations, then
\begin{itemize}
\item either $\Phi(M)$ is a cylinder,
\item or $\Iso_\Psi\FD=\FT\times Q$,
\item or $\Iso_\Psi\FD=\FT\times Q\times R$,
\end{itemize}
where $\times$ denotes the product of sets, 
$Q$ is a, possibly trivial, discrete group of translations, not
contained in $\FT$, and $R$ is the group generated by the
$180^\circ$-rotation around a fixed point, $z\rightarrow2z_0-z$, $z_0\in\Bcc$.
}

\Proof
1. W.l.o.g.\ we can assume that $\FR$ is the set of rotations around
the origin in $\Bcc$. If $\Iso_\Psi\FD\neq\FR$, then 
there exists $g\in\Iso_\Psi\FD$ s.t.\ $g(z)=az+b$ with
$b\neq0$ and, by Proposition~\ref{Cautom}, $|a|=1$. 
If $\tg(z)=az$, then
$\tg\in\FR$ and $(g\circ\tg\inv)(z)=z+b$, i.e., $\Iso_\Psi\FD$ contains
a pure translation. Lemma~\ref{Cautom} then
shows that $\Phi(M)$ is a cylinder.

2. Assume, that $\Iso_\Psi\FD$ contains a one-parameter group $\FT$ of
translations in $\Bcc$, i.e., 
\BEQ
\FT=\{g_r|g_r(z)=z+rv,\;r\in\Brr,v\in\Bcc\setminus\{0\}\},
\EEQ
and an automorphism $g(z)=az+b$ with $a\neq1$. Then
$|a|=1$ by Proposition~\ref{Cautom}.
If $b=0$, then $g$ is a rotation around the origin.
If $b\neq0$ we can, as in Proposition~\ref{Cautom},
apply the coordinate transformation $\tz=\frac{a-1}{b}z+1$, s.t.\ 
$g$ becomes a rotation around the origin, $g(\tz)=a\tz$.
In the new coordinates $\FT$ is still a
one-parameter group of pure translations. Its direction in the
$\tz$-plane is given by the vector $\tv=\frac{a-1}{b}v$.
Therefore, we can restrict
our investigation (after an additional rotation around the origin) to the case
that $\FT$ is the group of translations along the real axis and
$g$ is a rotation around the origin in $\Bcc$, $g(z)=az$.
Lemma~\ref{Cautom} shows, that then either $a=-1$ or
$\Phi(M)$ is a cylinder.
Therefore, if $\Phi(M)$ is not a cylinder, $Iso_\Psi\FD=T\times R$
or $\Iso_\Psi\FD=T$,
where $T$ is a group of translations and $\FT\subset T$. By the
same argument as in Lemma~\ref{Cautom}, for a noncylindrical surface, 
the orbit of $z=0$ under $\Iso_\Psi\FD$ cannot be the whole complex plane. This
shows, that $T=\FT\times Q$, where $Q$ is a, possibly trivial,
discrete group of translations.
\QED

\mycorollary{} {\em
Let $(M,\Phi)$, $(\Bcc,\Psi)$ and $\Iso_\Psi\FD$ be defined 
as in Theorem~\ref{Cautom}. Let
$E\diffz^2$ be the Hopf differential of $(\Bcc,\Psi)$. 
\begin{description}
\item[1.] If $\Iso_\Psi\Bcc$
contains the group of all
rotations around a fixed point, then, up to a biholomorphic change of
coordinates, we have
$E=d(z-z_0)^m$, where $d\in\Bcc\setminus0$ and $m=0,1,2,\ldots$ is an integer.
\item[2.] If $\Iso_\Psi\Bcc$ contains a one-parameter group $\FT$ of
translations, then, up to a biholomorphic change of coordinates, we
have $E\equiv1$.
\end{description}
}

\Proof 1. W.l.o.g.\ we can choose $z_0=0$. Let 
$\FR=\{g_\varphi\in\Aut\Bcc|g_\varphi(z)=e^{i\varphi}z,\;\varphi\in[0,2\pi)\}$
be the one-parameter group of rotations around the origin. 
By Corollary~\ref{associatedfamily}.1, $|E|$ is invariant
under all automorphisms in $\FR\subset\Iso_\Psi\FD$. Therefore, for each $\varphi\in\Brr$, we get
$E(g_\varphi(z))=e^{i\theta}E(z)$, where
$\theta=\theta(\varphi)\in\Brr$ depends linearly on $\varphi$.
It follows, since $E$ is holomorphic, that $\theta=m\phi$, $m$ a
nonnegative integer,
and, up to a biholomorphic change of coordinates, $E=dz^m$,
$d\in\Bcc$. Since by assumption $E\not\equiv0$, we have $d\neq0$.

2. By Corollary~\ref{associatedfamily}, $|E|$ and therefore also the
set of zeroes of $E$, is invariant under the group $\FT$. Therefore,
since the set of zeroes of a holomorphic function is discrete, $E$
cannot have any zeroes. It can therefore, by a
biholomorphic change of coordinates $\diff w^2=E(z)\diffz^2$, 
be transformed into $E\equiv1$.
\QED

\myremark{} The immersions considered in the theorem and in the
corollary will be
investigated in more detail in Section~\ref{SmythDel2}.

\newsection \label{uniqueness}
In the next sections we will investigate some properties of
the groups defined in Definition~\ref{secfund}. We begin with the following

\mylemma{} {\em 

a) Let $g\in\Aut_\pi\FD$ and $\hg\in\Aut M$ be as in \bref{eqa},
then $\hg$ is uniquely defined.

b) Let $\hg\in\Aut_\Phi M$ and $\tT\in\Aut\Psi(\FD)$ be as in \bref{eqb}, then
$\tT$ is uniquely defined.

c) Let $g\in\Aut_\Psi\FD$ and $\tT\in\Aut\Psi(\FD)$ be as in \bref{eqc}, then
$\tT$ is uniquely defined.
}

\Proof 
a) Assume $\hg$ and $\hg^\prime$ both satisfy \bref{eqa}, then
$\hg(\pi(z))=\hg^\prime(\pi(z))$ for all $z\in\FD$. This implies
$\hg=\hg^\prime$, since $\pi$ is surjective.

b) A proper Euclidean motion in $\threespace$ is
determined uniquely by its restriction to an affine two dimensional
subspace.  If we choose a point $z\in M$, then for each point $p$
of the affine tangent plane $\Phi(z)+\diff\Phi(T_zM)$,
$p=\Phi(z)+\diff\Phi(v)$, we have
\BEQ
\tT(p)=\Phi(\hg(z))+(\hg_\ast\diff \Phi)(v).
\EEQ
Therefore $\tT$ is uniquely determined.

c) Similarly.
\QED

\myremark{} 
It actually follows from the proof, that $\tT$ is already
determined by the restriction of $\hg$ to an arbitrary open 
subset of $M$.

\newsection \label{mapdefs}
Using Lemma~\ref{uniqueness} we define the following maps:

\mydefinition{}
\BEA
\bpi: \Aut_\pi\FD & \longrightarrow & \Aut M,\nonumber\\
\bpi: g & \longmapsto & \hg,
\EEA
where $g$ and $\hg$ are as in \bref{eqa},
\BEA
\phi: \Aut_\Phi M & \longrightarrow & \OAff(\threespace),\nonumber\\
\phi: \hg & \longmapsto & \tT,
\EEA
where $\hg$ and $\tT$ are as in \bref{eqb}, and
\BEA
\psi: \Aut_\Psi\FD & \longrightarrow & \OAff(\threespace),\nonumber\\
\psi: g & \longmapsto & \tT,
\EEA
where $g$ and $\tT$ are as in \bref{eqc}.

\separate By Lemma~\ref{uniqueness}, $\phi$ and $\psi$ are group homomorphisms.
Also note, that the images of $\phi$ and $\psi$ are contained in
$\Aut\Psi(\FD)$.

\mytheorem{} {\em 

a) The groups $\Aut_\Psi\FD$ and $\Aut_\Phi M$
are closed Lie subgroups of $\Aut\FD$ and $\Aut M$, respectively.

b) The maps $\bpi$, $\phi$ and $\psi$ are analytic homomorphisms of
Lie groups.
}

\Proof 
a) Let $g_n\in\Aut_\Psi\FD$ be a
sequence which converges to $g\in\Aut\FD$. Then $g_n$ converges
uniformly on each compact subset of $\FD$. In particular,
$\Psi\circ g_n=\tT_n\circ\Psi$ converges uniformly to $\Psi\circ
g$ on each sufficiently small closed ball around any point $z\in\FD$.
Therefore, also the differentials converge, whence
$(\tT_n)_\ast\diff_z\Psi=R_{\tT_n}\diff_z\Psi$ converges, where we have
written $\tT_n=(R_{\tT_n},t_{\tT_n})$ as in \bref{congdecomp}. 
This implies, that
$R_{\tT_n}$ converges to a rotation $R$ in $\threespace$. Since also
$\tT_n\circ\Psi=t_{\tT_n}+R_{\tT_n}\circ\Psi$ converges,
$t_{\tT_n}\rightarrow t$ for some $t\in\threespace$.
Altogether, this shows $\tT_n\rightarrow\tT=(R,t)$.
But now $\Psi\circ g_n\rightarrow\Psi\circ g=\tT\circ\Psi$. This
shows, that $\tT\in\Aut\Psi(\FD)$ and $g\in\Aut_\Psi\FD$. The argument
for $\Aut_\Phi M$ is similar.

b) We know, that $\Aut M$ and $\Aut\FD$ are Lie groups and, by the
argument above, we know that $\Aut_\Psi\FD$ and $\Aut_\Phi M$ are
closed subgroups of $\Aut\FD$ and $\Aut M$, respectively. Therefore,
with the induced topology $\Aut_\Psi\FD$ and $\Aut_\Phi M$ are Lie
groups. We show, that in this topology the maps $\phi$ and
$\psi$ are continuous, from which analyticity follows \cite[Th.~II.2.6]{He:1}.

Assume $g_n\rightarrow g$ in $\Aut_\Psi\FD$. Then $g_n$ converges to
$g$ uniformly on each compact subset of $\FD$. In particular,
$\Psi\circ g_n=\tT_n\circ\Psi$ converges uniformly to $\Psi\circ
g=\tT\circ\Psi$ on each sufficiently small closed ball around any
point $z\in\FD$.  By the proof of Lemma~\ref{uniqueness}, $\tT$ and
$\tT_n$ are uniquely determined by the restriction of $g$ and $g_n$ to
an arbitrary open subset of $\FD$. Therefore, $\tT_n$ converges to $\tT$. This
shows, that $\psi$ is continuous. For $\phi$ we proceed
analogously. For $\bpi$ the claim is trivial.
\QED

\newsection \label{immersedautom}
We will need the following

\mytheorem{} {\em 
We retain the notation of Section~\ref{secfund}.
If $M$ with the metric induced by $\Phi$ is complete, then
for every Euclidean motion $\tT\in\Aut\Psi(\FD)$ there exists 
a $g\in\Aut_\Psi\FD$, s.t.\ $\Psi\circ g=\tT\circ\Psi$, i.e.,
$\psi$ maps $\Aut_\Psi\FD$ onto $\Aut\Psi(\FD)$. The automorphism $g$
is unique up to multiplication with an element of $\Ker\psi$.
} 

\Proof 
Let $\tT$ be a Euclidean motion, which leaves the image $\Psi(\FD)$
invariant. Let us choose two arbitrary points $z_0$ and $z_1$ in
$\FD$, s.t.\ $\Psi(z_1)=\tT\Psi(z_0)$.
Then, since $\Psi$ is locally injective and conformal,
for each such pair $(z_0,z_1)$ there exists an open neighbourhood 
$U_0$ of $z_0$ and
an open neighbourhood $U_1$ of $z_1$, s.t.\
$\tT\circ\Psi(z)=\Psi(h(z))$, $z\in U_0$, defines an
orientation preserving isometry $h:U_0\rightarrow U_1$.

Since $M$ is complete, also $\FD$ is complete \cite[Prop.~I.10.6]{He:1}.
Thus, $\FD$ is an analytic, complete, simply connected manifold.
Therefore, by \cite[Sect.~I.11]{He:1},
the local isometry $h$ can be
extended to a unique, global, orientation preserving self-isometry 
$g\in\Aut\FD$ of $\FD$, s.t.\ $g|_{U_0}=h$. By the definition of
$h$ we have $\Psi(g(z))=\tT\Psi(z)$ on $U_0$. Since all occuring
maps are analytic, we get $\tT\circ\Psi=\Psi\circ g$ on $\FD$.
Uniqueness of $g$ up to an element of $\Ker\psi$ is trivial.
\QED

\separate The only group still to be discussed is $\Aut\Psi(\FD)$.
Unfortunately, in general $\Aut\Psi(\FD)$ doesn't seem to be closed in
$\OAff(\threespace)$.
In Section~\ref{regularsurfaces} we will give a simple condition on
$(M,\Phi)$, under which $\Aut\Psi(\FD)$ can be shown to be closed.

However, here we are able to conclude:

\mycorollary{} {\em 
If $(M,\Phi)$ is complete, then
$\psi:\Aut_\Psi\FD\rightarrow\Aut\Psi(\FD)$ is surjective. In
particular, $\Aut\Psi(\FD)\cong\Aut_\Psi\FD/\Ker\psi$ is a Lie group.
}

Theorem~\ref{immersedautom} and
Corollary~\ref{immersedautom} have well known equivalents for the map
$\pi$. The arguments leading to Eq.~\bref{AutMcong} prove the following

\myprop{} {\em Let $M$ be a Riemann surface with simply connected
cover $\pi:\FD\rightarrow M$. With the notation as above we have:

a) For every $\hg\in\Aut M$ there exists a
$g\in\Aut_\pi\FD$, s.t.\ $\pi\circ g=\hg\circ\pi$.

b) The map $\bpi:\Aut_\pi\FD\rightarrow\Aut M$ is surjective.
}

\newsection \label{groups}
We want to investigate, how the groups defined in
Section~\ref{secfund} are related to each other by the maps $\bpi$,
$\phi$ and $\psi$.

For every $g\in\bpi\inv(\Aut_\Phi M)$ we have
\BEQ \label{id1}
\Psi\circ g=\Phi\circ\pi\circ
g=\Phi\circ\bpi(g)\circ\pi=\phi(\bpi(g))\circ\Psi,
\EEQ
hence
\BEQ \label{inclimp}
\bpi\inv(\Aut_\Phi M)\subset\Aut_\Psi\FD
\EEQ
and $\psi=\phi\circ\bpi$ on $\bpi\inv(\Aut_\Phi M)$.

For $g\in\Ker\bpi$ we have $\Psi\circ g=\Phi\circ\pi\circ
g=\Phi\circ\pi=\Psi$. Therefore,
\BEQ
\Ker\bpi\subset\Ker\psi.
\EEQ
We also recall, that $\Ker\bpi=\Gamma$, the Fuchsian group of $M$. 

\mylemma{} {\em
Let $(M,\Phi)$ be a CMC-immersion with
$\Ker\psi=\Ker\bpi$. Then the following holds:
\begin{description}
\item[a)] $(\bpi)\inv(\Aut_\Phi M)=\Aut_\Psi\FD$.
\item[b)] $\phi:\Aut_\Phi M\longrightarrow\Aut\Psi(\FD)$ is an
injective group homomorphism.
\end{description}
If, in addition, $(M,\Phi)$ is complete, then:
\begin{description}
\item[c)]
The action of $\tT\in\Aut\Psi(\FD)$ can be lifted to an action on $M$, i.e.,
$\phi$ is surjective.
\item[d)] $\phi:\Aut_\Phi M\rightarrow\Aut\Psi(\FD)$ is a group isomorphism.
\end{description}
}

\Proof
a) Since $\Ker\psi=\Ker\bpi$, we have that $\Aut_\Psi\FD$ is in the
normalizer of $\Ker\bpi=\Gamma$, whence
$\Aut_\Psi\FD\subset\Aut_\pi\FD$, by~\bref{AutpiFD}. Therefore, for
$g\in\Aut_\Psi\FD$ we have
\BEQ
\Phi\circ\bpi(g)\circ\pi=\Phi\circ\pi\circ g=\Psi\circ g=\psi(g)\circ\Psi
=\psi(g)\circ\Phi\circ\pi,
\EEQ
and $\bpi(g)\in\Aut_\Phi M$ follows, i.e.,
$\Aut_\Psi\FD\subset\bpi\inv(\Aut_\Phi M)$. 
Now a) follows from Eq.~\bref{inclimp}.

b) From Eq.~\bref{id1} and a) it follows, that
$\psi=\phi\circ\bpi$ on $\Aut_\Psi\FD$, 
which implies $\Ker\phi=\{\id\}$ and therefore b).

c) Let $\tT\in\Aut\Psi(\FD)$. Since $(M,\Phi)$ is complete, there
exists, by Theorem~\ref{immersedautom}, $g\in\Aut_\Psi\FD$, s.t.\
$\tT=\psi(g)$. From a) it follows, that there exists $\hg\in\Aut_\Phi M$, s.t.\
$\hg=\bpi(g)$, whence $\tT=\phi(\hg)$. The map $\phi$ is therefore surjective.

d) follows from b) and c).
\QED

\myprop{} {\em

a) $\Ker\psi$ is a discrete subgroup of $\Aut_\Psi\FD$ and acts freely
and discontinuously on $\FD$.

b) $M^\prime=\Ker\psi\setminus\FD$ is a Riemann surface.

c) Let $\pi^\prime:\FD\rightarrow M^\prime$ denote the natural
projection. Then there exists an immersion $\Phi^\prime$ of $M^\prime$
into $\threespace$, s.t.\ the following diagram commutes:
\BEQ \label{MMpdiagram}
\setlength{\unitlength}{0.00083300in}%
\begin{picture}(3184,1113)(3851,-3841)
\thinlines
\put(5701,-3661){\vector( 3, 1){900}}
\put(4201,-3211){\vector( 1, 0){2400}}
\put(4201,-3361){\vector( 3,-1){900}}
\put(5731,-2800){\vector( 3,-1){870}}
\put(4051,-3286){\makebox(0,0)[b]{$\FD$}}
\put(4201,-3090){\vector( 3, 1){900}}
\put(5401,-3840){\makebox(0,0)[b]{$M^\prime=\Ker\psi\setminus\FD$}}
\put(6301,-3650){\makebox(0,0)[b]{$\Phi^\prime$}}
\put(6751,-3286){\makebox(0,0)[b]{$\threespace$}}
\put(5750,-3190){\makebox(0,0)[b]{$\Psi$}}
\put(4735,-3511){\makebox(0,0)[b]{$\pi^\prime$}}
\put(4440,-2911){\makebox(0,0)[b]{$\pi$}}
\put(6301,-2911){\makebox(0,0)[b]{$\Phi$}}
\put(5401,-2780){\makebox(0,0)[b]{$M=\Ker\bpi\setminus\FD$}}
\put(5401,-2800){\vector(0,-1){750}}
\end{picture}
\EEQ

d) For the CMC-immersion $(M^\prime,\Phi^\prime)$ as above we define
$\bpi^\prime$ and $\phi^\prime$ as in Section~\ref{mapdefs}. Then
\BEQ \label{kerbpikerpsi}
\Ker\bpi^\prime=\Ker\psi.
\EEQ
}

\Proof a) Since 
$\psi:\Aut_\Psi\FD\rightarrow\OAff(\threespace)$ is a continuous homomorphism
of Lie groups, $\Ker\psi$ is, with the induced topology, a Lie
subgroup of $\Aut_\Psi\FD$.
Therefore, if $\Ker\psi$ were nondiscrete, it would contain a
one-parameter subgroup $\gamma(t)$. Hence $\Psi(\gamma(t).z)=\Psi(z)$
for all $z\in\FD$ and all $t\in\Brr$. This implies $\gamma(t).z=z$ for
all $z\in\FD$ and all $t\in\Brr$, whence $\gamma(t)=I$ for all $t$, a
contradiction.

Now let us assume, that $g\in\Ker\psi$ has a fixed point $z_0\in\FD$.
Then
\BEA\label{one}
\Psi(g(z)) & = & \Psi(z)\kern1cm\mbox{\rm for all $z\in\FD$},\\
g(z_0) & = & z_0.\label{gptwo}
\EEA
Taking into account the injectivity of the derivative of $\Psi$ one
gets by differentiating Eq.~\bref{one} at $z=z_0$,
\BEQ \label{gpone}
g^\prime(z_0)=1.
\EEQ
In the case of $\FD=\Bcc$ one has $g(z)=az+b$, with $a,b\in\Bcc$. It
follows from Eqs.~\bref{gpone} and \bref{gptwo}, that $g=\id$.

In the case of $\FD$ being the unit circle we can view $g$ as an isometry
w.r.t.\ the Bergmann metric on $\FD$. This together with
\cite[Lemma~I.11.2]{He:1} implies again $g=\id$.

It remains to be proved, that $\Ker\psi$ acts discontinuously, i.e.\
that there is a point $z_0\in\FD$, s.t.\ the orbit of $\Ker\psi$
through $z_0$ is discrete.
For the unit circle this follows from the discreteness of $\Ker\psi$
and \cite[Theorem IV.5.4]{FaKr:1}. For $\FD=\Bcc$ it is trivial,
since then, by the arguments above, 
$\Ker\psi$ is a discrete group of translations.

b) Since by a), $\Ker\psi$ is a Fuchsian or elementary group, it
follows, that $M=\Ker\psi\setminus\FD$ is a Riemann surface (see
e.g.~\cite[Section~IV.5]{FaKr:1}).

c) From the definition of $\Ker\psi$ and $M^\prime$ it is clear that
$\Psi$ factors through $M^\prime$. This defines an immersion
$\Phi^\prime:M^\prime\rightarrow\threespace$. If
$\pi^\prime:\FD\rightarrow M^\prime$ is the natural projection, then
$\Phi^\prime\circ\pi^\prime=\Psi$ and~\bref{MMpdiagram} follows.

d) is clear from the definition of $M^\prime$.
\QED

\separate The last lemma shows, that for our purposes it is actually enough to
restrict our attention to surfaces with
\BEQ \label{simplycovered}
\Ker\psi=\Ker\bpi.
\EEQ
For these surfaces the conclusions of Lemma~\ref{groups} hold.

\newsection \label{nondiscrete}
The following proposition shows, what it means for the symmetry group
$\Aut_\Psi\FD$ that $\Aut\Psi(\FD)$ is not discrete.

\myprop{} {\em
Let $(M,\Phi)$ be a CMC-immersion with simply connected cover $\FD$,
which is complete w.r.t.\ the induced metric and admits a one
parameter group of Euclidean motions $P\subset\Aut\Psi(\FD)$. Then
$\Aut_\Psi\FD$ also contains a one parameter group.
}

\Proof Let $P=\{\tT_x,x\in\Brr\}$ be a
one parameter subgroup of $\Aut\Psi(\FD)$, where $\Aut\Psi(\FD)$
carries the Lie group structure stated in Corollary~\ref{immersedautom}.
Let $A\subset\FD$ be an open subset such that $\Psi$ is injective on
$A$. Let $a\in A$ be arbitrary. Since $\tT_0\Psi(a)=\Psi(a)\in\Psi(A)$, there
exists some $\epsilon>0$ and an open subset $A_\epsilon\subset A$, s.t.\
$\tT_x\Psi(A_\epsilon)\subset\Psi(A)$ for all $|x|<\epsilon$.
Therefore, by Theorem~\ref{immersedautom}, there exists an
automorphism $g_x\in\Aut\FD$, $|x|<\epsilon$, satisfying $\Psi\circ
g_x=\tT_x\circ\Psi$, which is unique up to multiplication with an
element in $\Ker\psi$. In addition, it follows from the proof of
Theorem~\ref{immersedautom}, that we can choose $g_x$ s.t.\ 
\BEQ \label{uniquechoice}
g_x(A_\epsilon)\subset A.
\EEQ
Since $\Ker\psi$ is discrete, the
condition~\bref{uniquechoice} determines $g_x$ uniquely, if $A$ is
small enough. This shows, that $g_{x+y}=g_xg_y$ for all sufficiently small
$x,y\in\Brr$. For $x\rightarrow0$ we have
$\tT_x\rightarrow\tT_0=I$, therefore $\Psi\circ g_x\rightarrow\Psi$
uniformly on $\FD$. This shows, that $g_x$ converges to an element of
$\Ker\psi$. Since $\Ker\psi$ acts freely on $\FD$, 
Eq.~\bref{uniquechoice} implies $g_x\rightarrow g_0=I$ for
$x\rightarrow0$. This shows that $\omega:x\rightarrow g_x$ is a
continuous and thus analytic (see~\cite[Th.~II.2.6]{He:1}) 
homomorphism from some interval $(-\tepsilon,\tepsilon)$, $\tepsilon>0$, into
$\Aut_\Psi\FD$. For an arbitrary $x\in\Brr$ we write 
$x=m\frac{\tepsilon}{2}+r$,
where $r\in [0,\frac{\tepsilon}{2})$ and $m\in\Bii$ are uniquely
determined. The definition
\BEQ
g_x=g_r(g_{\frac{\tepsilon}{2}})^m\in\Aut_\Psi\FD.
\EEQ
extends $\omega$ to a one-parameter subgroup of $\Aut_\Psi\FD$, which
finishes the proof.
\QED

\newsection \label{regularsurfaces}
It remains to be investigated under which circumstances the existence of a
cluster point of $\Aut\Psi(\FD)$ implies $\dim\Aut\Psi(\FD)\geq1$.
This is certainly the case, if
$\Aut\Psi(\FD)$ is closed in $\OAff(\threespace)$.

To this end we introduce the notion of an admissible immersion.

\mydefinition{} 
Let $(M,\Phi)$ be an immersed manifold in $\threespace$. A point
$p\in\Phi(M)$ is called {\em admissible}, if
there is an open neighbourhood $U$ of $p$ in $\threespace$,
s.t.\ the intersection $\Phi(M)\cap U$ is closed
in $U$. The immersion $(M,\Phi)$ is called admissible, if
$\Phi(M)$ contains at least one admissible point.

\separate 
We think it is fair to say that, basically, every surface of interest
is admissible.  Most surfaces studied actually belong to the smaller
class of locally closed surfaces (see \cite[II.2]{La:1}), for which each
point of the image is admissible.
Among the locally closed surfaces are e.g.\ the 
immersed surfaces with closed image $\Phi(M)$ in $\threespace$, 
especially compact submanifolds of $\threespace$,
and immersed surfaces $(M,\Phi)$, for which $\Phi$ is proper
(see e.g.~\cite[I.2.30]{Sp:1}).

Also note that, geometrically speaking, a surface has to return
infinitely often to each neighbourhood of each of its nonadmissible points.
A nonadmissible surface is therefore in a sense a
two-dimensional analog of a Peano curve.

\myremark{} It is important to note, that admissibility is a property
of  the image $\Phi(M)$ of the immersion $\Phi$.
We don't claim that it is preserved under isometries (see
Section~\ref{associatedfamily}).
In particular, for an admissible surface it may well be, that 
not all members of the associated family are admissible. 

\separate
The definition of admissible surfaces allows us to describe a large
class of surfaces, for which $\Aut\Psi(\FD)$ is closed.

\mytheorem{} {\em If $(M,\Phi)$ is a complete, admissible surface in
$\threespace$ and $(\FD,\Psi)$
is the simply connected cover of $M$ with the covering immersion
$\Psi=\Phi\circ\pi$, then the group $\Aut\Psi(\FD)$ is closed in
$\OAff(\threespace)$.
}

\underline{Proof:}
Let $\tT_n\in\Aut\Psi(\FD)$ be a sequence of symmetry transformations
of $\Psi(D)$ which converges to $\tT\in\OAff(\threespace)$. Therefore, also
the sequence $\tT_n\inv$ converges in $\OAff(\threespace)$. 

Since $(M,\Phi)$ is admissible, there exists an admissible point 
$p\in\Psi(D)$, together with an open ball $B(p,\epsilon)$ of radius
$\epsilon<1$ around $p$ in $\threespace$, 
s.t.\ $B(p,\epsilon)\cap\Psi(\FD)$ is closed in $B(p,\epsilon)$.

W.l.o.g. we can assume that $p$ and the whole bounded sequence
$\{\tT_n\inv(p)\}$ lies in $B(0,{1\over2})$. Otherwise we first apply a scaling
transformation of $\threespace$, which changes neither the admissibility of
$(M,\Phi)$ nor the group structure of $\Aut\Psi(\FD)$.

We take $N\in\Bnn$ s.t.\ 
$\|\tT-\tT_n\|<{\epsilon\over3}$ for $n\geq N$,
where $\|\cdot\|$ denotes the operator norm.

We choose $p^\prime=\tT_N\inv(p)$ and $z^\prime\in\FD$, s.t.\
$p^\prime=\Psi(z^\prime)$. Since
$\tT_N$ is a Euclidean motion we have that
$\tT_N(B(p^\prime,\epsilon))=B(p,\epsilon)$.

For all $q\in B(p^\prime,{\epsilon\over3})\cap\Psi(\FD)$ we get
$|q|\leq |q-p^\prime|+|p^\prime|<{\epsilon\over3}+{1\over2}<1$,
therefore, if $n\geq N$,
\BEQ
|\tT_n(q)-p|\leq
|\tT_n(q)-\tT(q)|+|\tT(q)-\tT(p^\prime)|+|\tT(p^\prime)-\tT_N(p^\prime)|<
\frac{5}{6}\epsilon.
\EEQ
Thus we have $\tT_n(q)\in B(p,\frac{5}{6}\epsilon)\cap\Psi(\FD)$ for
all $n\geq N$ and $\overline{\{\tT_n(q),n\geq N\}}\subset B(p,\epsilon)$.
Now we know, that $B(p,\epsilon)\cap\Psi(\FD)$ is closed in
$B(p,\epsilon)$, and that $\tT_n(q)$ converges by assumption to
$\tT(q)\in\threespace$. Therefore, the limit $\tT(q)$ is in
$\Psi(\FD)$. Since $\Psi$ is an immersion, there exists an open
neighbourhood $A$ of $z^\prime$ which is mapped into
$B(p^\prime,\frac{\epsilon}{3})$ by $\Psi$.
Therefore, if we choose $z\in\FD$, s.t.\ $p=\Psi(z)$, then 
$\tT$ induces an isometry $g$ of $A$ onto an open
neighbourhood $B\subset\FD$ of $z$, s.t.\ $\tT\circ\Psi=\Psi\circ g$ on
$A$. By \cite[Sect.~I.11]{He:1}, this isometry can be extended globally to a
unique automorphism $g\in\Aut\FD$. Since all
maps are analytic and globally defined, the relation
$\tT\circ\Psi=\Psi\circ g$ holds on the whole of $\FD$.
It follows, that $g\in\Aut_\Psi\FD$ and $\tT\in\Aut\Psi(\FD)$.
\QED

\mycorollary{}
Let $(M,\Phi)$ and $\FD$ be as in Theorem~\ref{regularsurfaces}.
If $\Aut\Psi(\FD)$ is nondiscrete, then also
$\Aut_\Psi\FD$ is nondiscrete.

\Proof By the assumptions,
$\Aut\Psi(\FD)$ is closed in $\OAff(\threespace)$
and nondiscrete. It therefore contains a
one parameter group and the corollary follows from
Proposition~\ref{nondiscrete} above.\QED

Finally, we state the following result for the translational parts of
the elements of $\Aut\Psi(\FD)$:

\myprop{} {\em
If $(M,\Phi)$ is admissible and complete, and $\Aut\Psi(\FD)$ is
discrete, then the set
$\FL=\{t|\tT=(R,t)\in\Aut\Psi(\FD)\}$ 
of translations is discrete.
}

\Proof Assume, $\FL$ is not discrete. Then there
exists a $t\in\FL$ and a sequence $\{t_n\}\subset\FL$, $t_n\neq t$, s.t.\
$t_n$ converges to $t$. Since the set of rotations is compact, we obtain
a subsequence $\tT_n=(R_n,t_n)\in\Aut\Psi(\FD)$, which
converges in $\OAff(\threespace)$. 
Since by Theorem~\ref{regularsurfaces}, $\Aut\Psi(\FD)$ is closed, 
this sequence converges in $\Aut\Psi(\FD)$ to some $\tT$. But
then $\tT_n=\tT$ for sufficiently large $n$, since $\Aut\Psi(\FD)$ is
discrete. This shows, that $t_n=t$ for sufficiently large
$n$, a contradiction.
\QED

\newsection \label{DelundSmyth}
As examples for the discussion in this chapter,
let us investigate two well known classes of CMC-surfaces,
the Delaunay and the Smyth surfaces.

We recall that a Delaunay surface is defined as a complete, immersed surface of
constant mean curvature which is generated in $\threespace$ 
by rotating a curve around a given axis. 
Clearly, every sphere and every cylinder is a Delaunay surface.
We will restrict the definition to those surfaces which have mean
curvature $H=-\frac{1}{2}$ and exclude the degenerate case of the sphere.

Let us translate two well known facts about Delaunay surfaces
into our language:

\myprop{} {\em
1. Let $(M,\Phi)$ be a noncylindrical CMC-immersion with universal 
covering immersion $(\FD,\Psi)$, s.t.\ $\Phi(M)$ is a Delaunay
surface. Then $\Phi(M)$ is generated by rotating 
the roulette of an ellipse (unduloid) or a hyperbola (nodoid) 
along the line on which the conic rolled.

2. Let $S\subset\threespace$ be a Delaunay surface. Then there exists
a CMC-immersion $(M,\Phi)$ with universal covering immersion $(\FD,\Psi)$,
$\FD=\Bcc$, s.t.\
\begin{itemize}
\item $S=\Phi(M)$,
\item $\Aut_\Psi\FD$ contains a one parameter group $\FT$ of
translations, wich is mapped by the surjective homomorphism 
$\psi:\Aut_\Psi\FD\rightarrow\Aut\Psi(\FD)$ to the group of 
rotations around the axis of revolution of the Delaunay surface.
\end{itemize}
}

\Proof 1. well known, see e.g.~\cite{Ee:1}.

2. such immersions (with $\Phi=\Psi$, $M=\FD=\Bcc$)
are explicitly constructed in~\cite{Sm:2}.
\QED

\separate
For later use, we collect some definitions for surfaces of
revolution (see e.g.~\cite{Sto:1}): 
If $S$ is a surface of revolution, generated by rotating
the plane curve $C\subset\Bcc$ around the axis $A$, then $A$ is called the
{\em axis of revolution}, the images of $C$ under a rotation around
$A$ are called {\em meridians}, and the intersection circles of $S$ with any
plane perpendicular to $A$ are called the {\em parallels} of $S$. 

\myremark{}
1. While it is clear for the unduloid, that the roulette of an ellipse
gives a periodic curve along $A$, this is not clear for the
nodoid. Here, in order to get a complete noncompact surface w.r.t.\
the induced metric, one has to continue the roulette periodically
along $A$. This is possible in a unique way (see~\cite{Ee:1}) and
gives a periodic curve along $A$. 

2. Besides of being periodic, the roulette of an ellipse or hyperbola
has another important property: In each of its periods there is a
unique point of maximal distance from the line $A$ on which the conic
rolls. The roulette is symmetric w.r.t.\ the reflection at any line 
perpendicular to $A$, which passes through such a point of maximal
distance from $A$. The Delaunay surface is therefore invariant under a
$180^\circ$-rotation around any axis which is perpendicular to the axis
of revolution $A$ and passes through a parallel of maximal radius.

\separate
Let us recall again, that in the discussion of Delaunay surfaces
we exclude the catenoid, the
roulette of the parabola, which is a minimal surface.
Since we also exclude the sphere, we get the following result, the
proof of which is an exercise in elementary geometry:

\mylemma{} {\em 
Each Delaunay surface determines its axis of revolution uniquely.
}

\separate
Proposition~\ref{DelundSmyth} gives the following

\mycorollary{} {\em
Let $(M,\Phi)$ be as in Proposition~\ref{DelundSmyth}. Let $A$ be the
generating axis of the Delaunay surface $\Phi(M)$. Then, as a set,
$\Aut\Psi(\FD)$ can be written as
\BEQ
\Aut\Psi(\FD)=\FR\times\tQ\times\{I,\tR\},
\EEQ
where $\FR$ is the one-parameter group of rotations around $A$,
$\tQ$ is a nontrivial discrete group of translations along $A$,
and $\tR$ is a $180^\circ$-rotation around an axis which is
perpendicular to $A$.
}

\Proof
By the remark above, every Delaunay surface is generated by a periodic
function that has only one maximum in every period interval.
Let $\FP\subset\threespace$ 
be the union of all parallels of $\Phi(M)$ which have
maximal radius. Denote by $I(\FP)$ the set of all proper Euclidean
motions, which leave $\FP$ invariant.
$\FP$ consists of disjoint circles which all lie on the same cylinder
around $A$. Let $\FC\subset A$ denote the set of centers of these
circles. Then, clearly, $I(\FP)$ is the set of all elements of
$\OAff(\threespace)$, which leave $A$ and $\FC$ invariant.
By the second part of Remark~\ref{DelundSmyth}, 
the map $\FC$ consists of 
equidistant points along $A$. It follows, that
$I(\FP)$ is generated by
\begin{itemize}
\item the one-parameter group $\FR$ of rotations around $A$,
\item a discrete group $\tQ$ of translations along $A$, which is given
by $\FC$, i.e., the periodicity of the roulette,
\item the $180^\circ$-rotation $\tR$ around an axis, which passes
through an arbitrary fixed point on $\FP$ and the center of the
corresponding parallel.
\end{itemize}
We therefore have
\BEQ
I(\FP)=\FR\times\tQ\times\{I,\tR\},
\EEQ
where `$\times$' denotes the product of sets.
By Lemma~\ref{DelundSmyth}, a Delaunay surface is left invariant
precisely by those proper Euclidean motions,
which preserve the axis of revolution and
map a meridian into another meridian. In particular, if
$p\in\Phi(M)$ and $\tT\in\Aut\Psi(\FD)$, then $p$ and $\tT(p)$ have
the same distance from $A$.
The set $\FP$ is the set of all
points on $\Phi(M)$, which have maximal distance from the
axis $A$. Therefore, every element of
$\Aut\Psi(\FD)$ leaves also the set $\FP$ invariant, i.e.,
$\Aut\Psi(\FD)\subset I(\FP)$.
As a surface of revolution around $A$, the Delaunay surface is
certainly invariant under the group $\FR$ of rotations around $A$.
In addition, by the second part of Remark~\ref{DelundSmyth}, 
for a Delaunay surface,
the meridians are periodic along the axis of revolution and
symmetric w.r.t.\ a $180^\circ$-rotation around any axis which is
perpendicular to $A$ and intersects $\FP$.
Therefore, $\Phi(M)=\Psi(\FD)$ is invariant also under the group $\tQ$ and the
rotation $\tR$ defined above.
This gives $I(\FP)\subset\Aut\Psi(\FD)$, which finishes the proof.
\QED

\separate
Smyth \cite{Sm:2} introduced for every integer $m\geq0$ a
one-parameter family of conformal immersions
\BEQ
\Psi^m_c\colon \Bcc\longrightarrow\threespace,\kern1cm c\in\Bcc\setminus\{0\},
\EEQ
with constant mean curvature, s.t.\ the induced metric is complete
and invariant under the one-parameter group of rotations around $z=0$
in $\Bcc$. We
will call these surfaces Smyth surfaces. The Hopf differential of
$(\Bcc,\Psi^m_c)$ is $cz^m\diffz^2$, and therefore each
$(\Bcc,\Psi^m_c)$ has an umbilic of order $m$ at the origin.
For $m=0$ the family $\Psi^m_c$ contains the cylinder.
We will call a Smyth surface {\em nondegenerate},
if its image in $\threespace$ is not a cylinder.

Two surfaces in $\threespace$ will be called congruent if they are related
by a proper Euclidean motion of $\threespace$.
Recall also the definition of the associated family in
Section~\ref{associatedfamily}.

The following results were proved by Smyth~\cite{Sm:2}:

\mytheorem{} {\em
Let $(M,\Phi)$, with covering immersion $(\FD,\Psi)$,
be a complete, immersed surface of constant mean curvature,
admitting a one-parameter group of self-isometries. 
Then the following holds:
\begin{description}
\item[1.] The simply connected cover of the Riemann surface $M$ 
is $\FD=\Bcc$.
\item[2.] The associated family of $(\FD,\Psi)$ contains either a
Delaunay or a Smyth surface, i.e., $(\FD,\Psi)$ is isometric to the
simply connected cover of a Delaunay or a Smyth surface.
\item[3.] The surface $(\FD,\Psi)$ admits a one-parameter group $P$ of
self-isometries which is
\begin{description}
\item[a)] a one-parameter group of translations ($P\cong\Brr$) 
in case of the Delaunay surfaces,
\item[b)] a one-parameter group of rotations around a fixed point in $\Bcc$
($P\cong S^1$) in case of the Smyth surfaces.
\end{description}
\end{description}
}

\newsection \label{SmythDel2} 
With the results in the previous section 
we can also easily derive the uniformization 
of Delaunay and Smyth surfaces.

\myprop{} {\em
1. Each Delaunay surface is conformally equivalent to the cylinder, i.e.,
its simply connected cover is $\Bcc$ and the
Fuchsian group is a one-parameter group of translations.

2. Each nondegenerate Smyth surface
is conformally equivalent to $\Bcc$. 

3. For a nondegenerate Smyth surface we have
\BEQ
\Ker\bpi=\Gamma=\{\id\}.
\EEQ
}

\Proof
We already know from Theorem~\ref{DelundSmyth}, that both, Delaunay
and Smyth surfaces, have, as Riemann surfaces, the 
simply connected cover $\Bcc$. Therefore, by Lemma~\ref{Riemannian},
they are biholomorphically equivalent to the plane, the cylinder or
a torus and the Fuchsian group is either trivial or it is a discrete
group of translations.

1. Delaunay surfaces are surfaces of revolution, i.e., there
exists a conformal immersion $\Psi:\Bcc\rightarrow\threespace$, which
is invariant under a discrete one-dimensional lattice $\Gamma$ of translations
in $\Bcc$, s.t.\ $\Psi(\FD)$ is the Delaunay surface. 
The Fuchsian group of the
underlying Riemann surface contains therefore a group of
translations. Since Delaunay surfaces are noncompact, they are
biholomorphically equivalent to the cylinder.

2. Let $\Psi:\Bcc\rightarrow\threespace$ be an immersion
s.t.\ $\Psi(\FD)=\Psi_c^m(\FD)$ is a Smyth surface for some parameters
$m\in\Bnn$, $c\in\Bcc$. Assume, that the Fuchsian group of the surface
contains a nontrivial translation. Then, by Theorem~\ref{DelundSmyth},
the group $\Aut_\Psi\FD$ satisfies the assumptions of the first part
of Lemma~\ref{Cautom}. Therefore, $\Psi(\FD)$ is a cylinder.
For a nondegenerate Smyth surface this gives
$M=\Bcc=\FD$, i.e., the Smyth surface is conformally equivalent to the
complex plane.

3. From 2.\ it follows, that $\pi=\id$ and $\Aut_\pi\FD=\Aut\FD=\Aut
M$, therefore $\Ker\bpi=\Gamma=\{\id\}$.
\QED

\separate
We continue with the following

\mylemma{} {\em
1. For each noncylindrical Delaunay surface $S\subset\threespace$,
there exists a CMC-immersion $(M,\Phi)$ with $\Phi(M)=S$, s.t.\
\BEQ \label{beh1}
\Ker\psi=\Ker\bpi,
\EEQ
\BEQ \label{beh2}
\Aut_\Phi M\cong\Aut\Psi(\FD),
\EEQ
and, as a product of sets, we have
\BEQ \label{beh3}
\Aut_\Psi D=\Iso_\Psi D=\FT\times Q\times R,
\EEQ
where $\FT\subset\Aut\Bcc$ is a one-parameter group of translations,
$Q\subset\Aut\Bcc$ is a discrete group of translations with one generator, and
$R=\{I,R_\pi\}\subset\Aut\Bcc$ is the group generated by the inversion
$R_\pi:z\mapsto-z$.

2. If $(\Bcc,\Psi=\Psi_c^m)$ is a nondegenerate Smyth surface, then
\BEQ
\Ker\psi=\Ker\bpi=\{\id\}
\EEQ
and
\BEQ \label{SmythAuts}
\Aut\Psi(\FD)\cong\Aut_\Psi\FD=\FR,
\EEQ
where $\FR$ is a finite group of rotations around $z=0$ in $\Bcc$.
}

\Proof
1. For a Delaunay surface, there exists, by Proposition~\ref{DelundSmyth},
a universal covering immersion
$\Psi:\Bcc\rightarrow\threespace$, s.t.\ the set
$\Aut_\Psi\FD\subset\Iso_\Psi\FD$
contains a one-parameter group $\FT\cong\Brr$ of translations.
The group $\FT$ is then mapped by the surjective homomorphism
$\psi:\Aut_\Psi\FD\rightarrow\Aut\Psi(\FD)$, defined in
Section~\ref{groups}, into the set of rotations
$\FR$ around $A$, the axis of revolution. By Proposition~\ref{groups} we can
further assume, that Eq.~\bref{beh1} holds. 
Therefore, $\FT$ contains the discrete Fuchsian group $\Gamma$ of the
Delaunay surface. By a rotation of the coordinate system in $\Bcc$, we
can choose
$\FT$ as the group of translations along the imaginary axis in $\Bcc$.
If we identify $M=\Gamma\setminus\Bcc$ with a fundamental domain
of $\Gamma$ in $\Bcc$, then $\Psi$ maps this region, a strip parallel
to the real axis, conformally to the Delaunay surface. The meridians
of $\Phi(M)$ are identified with lines in $M$ which are parallel to
the real axis, the parallels of $\Phi(M)$ are identified with the
intersection of $M$ with lines parallel to the imaginary
axis. In particular, the set $\FP$, which was defined in the proof of
Corollary~\ref{DelundSmyth}, is identified with the intersection of
$M$ with an
equidistant set $\FL$ of lines parallel to the imaginary axis. By a
translation of the coordinate system in $\Bcc$ we can, in addition, assume
$0\in\FL$.

By Eq.~\bref{beh1} and Lemma~\ref{groups}, we have that $\phi:\Aut_\Phi
M\rightarrow\Aut\Psi(\FD)$, defined in Section~\ref{groups}, is an
isomorphism of Lie groups. Therefore, Eq.~\bref{beh2} holds.
With the definitions in Corollary~\ref{DelundSmyth}
we can descibe the group $\Aut_\Phi M$, and, by the identification of
$M$ with a subset of $\Bcc$, also $\Aut_\Psi\FD$, explicitly:
\begin{itemize}
\item The group $\FR\subset\Aut\Psi(\FD)$ is identified in $\Aut M$ with 
the set of translations parallel to the imaginary axis in $\Bcc$. We
therefore have $\FT=\psi\inv(\FR)$.
\item The group $\tQ\subset\Aut\Psi(\FD)$ is identified in $\Aut M$ with the
group $Q$ of translations in $\Bcc$ which leave $\FL$ invariant. We
therefore have $\psi\inv(\tQ)=Q$.
\item The rotation $\tR$ is identified in $\Aut M$ 
with a $180^\circ$-rotation
around an arbitrary fixed point $z\in\FL$. We choose $z=0$. Then $\tR$
is identified with $R_\pi:z\rightarrow-z$ up to an automorphism in
$\Gamma\subset\FT$.
We therefore have $\psi\inv(\FR\times\{I,\tR\})=\FT\times\{I,R_\pi\}$.
\end{itemize}
Since $\psi$ is a surjective homomorphism and since
$\Aut\Psi(\FD)=\FR\times\tQ\times\{I,\tR\}$, we have
\BEQ
\Aut_\Psi\FD=\psi\inv(\Aut\Psi(\FD))=\FT\times Q\times\{I,R_\pi\}.
\EEQ
Since with Theorem~\ref{Cautom},
\BEQ
\FT\times Q\times\{I,R_\pi\}\subset\Aut_\Psi\FD\subset\Iso_\Psi\FD\subset
\FT\times Q\times\{I,R_\pi\},
\EEQ
we get $\Iso_\Psi\FD=\FT\times Q\times\{I,R_\pi\}$,
and therefore Eq.~\bref{beh3}.

2. For the immersions $\Psi_c^m$,
the metric $\frac{1}{2}e^u\diffz\diffzbar$ and, by
Eq.~\bref{Gausscurvature}, also $|E|$ is invariant under the
one-parameter group of rotations around $z=0$.
Therefore, by Corollary~\ref{secfund} and Theorem~\ref{Cautom}, 
we have that either $\Aut_\Psi\FD$ is contained in the group of 
rotations around a fixed point, or $\Psi(\FD)$ is a cylinder.
Therefore, for a nondegenerate Smyth surface,
$\Ker\psi\subset\Aut_\Psi\FD$ consists only of
rotations. But since, by Proposition~\ref{groups}, 
$\Ker\psi$ acts freely on $M=\Bcc$, this implies
that $\Ker\psi=\{\id\}=\Ker\bpi$ and, with Corollary~\ref{immersedautom},
$\psi:\Aut_\Psi\FD\rightarrow\Aut\Psi(\FD)$ is an isomorphism of Lie
groups, i.e.,
\BEQ
\Aut\Psi(\FD)\cong\Aut_\Psi\FD.
\EEQ
Let $g\in\Aut_\Psi\FD$. Then  $g$ is a rotation around $z=0$, i.e.,
$g(z)=e^{i\theta}z$ for some $\theta\in[0,2\pi)$. By Eq.~\bref{Etrafo1} and
Corollary~\ref{Cautom}, we get that $e^{i(m+2)\theta}=1$ for some integer
$m\geq0$. This shows, that $\Aut_\Psi\FD$ is a discrete group of rotations
around $z=0$ in $\Bcc$, finishing the proof.
\QED

\separate In Section~\ref{Smythexample} 
we will revisit the Smyth surfaces.

\section{Automorphisms in the DPW approach}
\label{DPWautom}\message{[DPWautom]}
We are working in the framework of \cite{DoPeWu:1} and \cite{DoHa:1}.
For notational conventions see the appendix of \cite{DoHa:1}.

In this chapter we consider only complete CMC-immersions $(M,\Phi)$ which
satisfy $\Ker\psi=\Ker\bpi$.  Therefore, by
Lemma~\ref{groups}, every $\tT\in\Aut\Psi(\FD)$ is associated with a
$\hg\in\Aut_\Phi M$. We have the following commutative
diagram of surjective group homomorphisms

\begin{center}
\unitlength0.5cm
\begin{picture}(3.5,3)
\put(-6.3,2.8){$\Aut_\Psi\FD=(\bpi\inv)(\Aut_\Phi M)$}
\put(3.2,2.8){$\Aut\Psi(\FD)$}
\put(-0.3,0){$\Aut_\Phi M$}
\put(0.4,1.5){$\bpi$}
\put(2,3.2){$\psi$}
\put(1.6,3){\vector(1,0){1.4}}
\put(0.2,2.5){\vector(0,-1){1.5}}
\put(1.5,1){\vector(1,1){1.5}}
\put(2.5,1.3){$\phi$}
\end{picture}
\end{center}

and $\phi$ is a group isomorphism.

\newsection \label{33}
Let $\hN$ be the Gau\ss\ map of $M$. 
We lift $\hN$ to $\FD$:

\begin{center}
\unitlength0.5cm
\begin{picture}(3.5,3)
\put(0,2.8){$\FD$}
\put(3.2,2.8){$S^2$}
\put(-1.2,0){$M=\Gamma\setminus\FD$}
\put(0.4,1.5){$\pi$}
\put(2,3.2){$N$}
\put(1,3){\vector(1,0){2}}
\put(0.2,2.5){\vector(0,-1){1.5}}
\put(1.5,1){\vector(1,1){1.5}}
\put(2.5,1.3){$\hN$}
\end{picture}
\end{center}

For $\gamma\in\Gamma$ we get $N\circ\gamma=N$ and for $g\in
\Aut_\Psi\FD$, $\tT=\psi(g)=(R_{\tT},t_{\tT})$, we get 
\BEQ \label{normtransform}
N\circ g=R_{\tT}\circ N.
\EEQ 
To translate everything into the language of Lie groups we note
$S^2\cong\LieSO(3)/\LieSO(2)\cong\LieSU(2)/\LieU(1)$. 
{}From now on we use, without changing the
notation, the spinor representation $J:\threespace\rightarrow\Liesu(2)$, of
$\threespace$, which is defined by
\BEQ
r\in\threespace\longmapsto
-\frac{i}{2}\vecr\cdot\vecsigma,\kern1cm\vecsigma=(\sigma_1,\sigma_2,\sigma_3),
\EEQ
$\sigma_i$, $i=1,2,3$, being the Pauli matrices.

Using the universal covering immersion
$\Psi:\FD\rightarrow\threespace$ we can lift the Gau\ss\ map
$N:\FD\rightarrow S^2$
to a frame $\hF=(e^{-\frac{u}{2}}\Psi_x,e^{-\frac{u}{2}}\Psi_y,N):
\FD\rightarrow\LieSO(3)$.
We will also normalize the frame, s.t.\ $\hF(0)=I$.

Consider now $g\in\Aut_\Psi\FD$. Then~\bref{normtransform} implies,
that also the tangent space at $\Psi(z,\zquer)$ is moved by $R_{\tT}$ to the
tangent space at $\Psi(g(z),\overline{g(z)})$. 
Since $R_{\tT}$ is orientation
preserving, the axes can be aligned by a rotation fixing $N$. Hence,
\BEQ \label{hFtrafo}
(\hF\circ g)(z,\zquer)=R_{\tT}\hF(z,\zquer)\hk(g,z,\zquer).
\EEQ
Note, that here $\hk$ is differentiable in $(z,\zquer)$. Since with
$g\mapsto\psi(g)$ also $g\mapsto R_{\tT}$ is a homomorphism, it is easy
to verify, that $\hk$ is a cocycle, i.e.,
\BEQ \label{hkcocycle}
\hk(g_2\circ
g_1,z,\zquer)=\hk(g_1,z,\zquer)\hk(g_2,g_1(z),\overline{g_1(z)}).
\EEQ
Since $\FD$ is simply connected, $\hF$ can be lifted to a frame $\tF$
taking values in $\LieSU(2)$. This lift is unique, if we require
\BEQ \label{tFinitial}
\tF(0)=I.
\EEQ
Taking a preimage $\tchi(g)\in\LieSU(2)$ of $R_{\tT}\in\LieSO(3)$ and
$k(g,z,\zquer)\in\LieU(1)$ of $\hk(g,z,\zquer)\in\LieSO(2)$ gives
\BEQ \label{frametrafo}
(\tF\circ g)(z,\zquer)=\tchi(g)\tF(z,\zquer)k(g,z,\zquer).
\EEQ
Clearly, $\tchi(g)$ and $k$ are determined up to a factor $\pm
I$. In particular, $k$ is differentiable in $(z,\zquer)$.
The homomorphism property of $R_{\tT}$ and the cocycle condition for
$\hk$ only lift to
\BEQ \label{tchihomom}
\tchi(g_2\circ g_1)=\pm\tchi(g_2)\tchi(g_1),
\EEQ
\BEQ \label{kcocycle}
k(g_2\circ g_1,z,\zquer)=\pm k(g_1,z,\zquer)k(g_2,g_1(z),\overline{g_1(z)}).
\EEQ
Although we cannot expect $\tchi(g)$ to be a homomorphism, at least
for translations in $\Aut_\Psi\FD$, $\FD=\Bcc$, we have

\mylemma{} {\em
Let $\Psi:\Bcc\rightarrow\threespace$ be a CMC-immersion. Define the
frame $\hF=(e^{-\frac{u}{2}}\Psi_x,e^{-\frac{u}{2}}\Psi_y,N):
\Bcc\rightarrow\LieSO(3)$ and its lift $\tF:\Bcc\rightarrow\LieSU(2)$ as
above. Let $T$ be the subgroup of translations in 
$\Aut_\Psi\Bcc$. Then $\hk(g,z,\zquer)\equiv I\in\LieSO(3)$ 
in Eq.~\bref{hFtrafo} for all $g\in T$. 
If one lifts $\hk(g,z,\zquer)$ to $k(g,z,\zquer)\equiv I\in\LieSU(2)$ for all
$g\in T$, then, after restriction to $T$, $\tchi(g)$ becomes a
homomorphism of groups from $T$ into $\Aut\Psi(\Bcc)$.
}

\Proof 
If $g\in T\subset\Aut_\Psi\Bcc$, then $\Psi\circ g=\tT\circ\Psi$ for
some $\tT=(R_{\tT},t_{\tT})\in\OAff(\threespace)$. 
{}From this it follows for the derivative of $\Psi$, since $g$ is a translation:
\BEQ
D(\Psi\circ g)(z,\zquer)=R_{\tT}(D\Psi(g(z),\overline{g(z)})\cdot
Dg(z,\zquer))=R_{\tT}(D\Psi(z,\zquer)).
\EEQ
Using Eq.~\bref{utrafo} and $g^\prime(z)=1$, we get
\BEQ
\hF\circ g=(e^{-\frac{u}{2}}R_{\tT}\Psi_x,
e^{-\frac{u}{2}}R_{\tT}\Psi_y,R_{\tT}N)=R_{\tT}\hF,
\EEQ
i.e., $\hk(g,z,\zquer)\equiv I$. From  this it follows, that we can
lift $\hk(g,z,\zquer)$ for all $g\in T$ to $k(g,z,\zquer)\equiv
I$. Since then, by Eq.~\bref{frametrafo},
\BEQ
(\tF\circ g)(z,\zquer)=\tchi(g)\tF(z,\zquer),
\EEQ
we have that $\tchi:T\rightarrow\Aut\Psi(\Bcc)$ is a homomorphism of groups.
\QED

\myremark{} 
In particular, if $M=\Gamma\setminus\FD$ is a torus or a Delaunay
surface, then the group $\Gamma$ consists of translations in
$\FD=\Bcc$. Lemma~\ref{33} then shows, that we can 
assume $k(\gamma,z,\zquer)=I$ for all
$\gamma\in\Gamma\subset\Aut_\Psi\Bcc$. Thus, $\tchi$ maps $\Gamma$
and the subgroup of all translations in $\Aut_\Psi\FD$ homomorphically
into $\Aut\Psi(\Bcc)$.

\newsection \label{Chistetig}
Since $S^2\cong\LieSU(2)/\LieU(1)$ is a compact symmetric space, we
have an associated Cartan decomposition $\Liesu(2)=\Liek+\Liep$, where
in a suitable matrix representation $\Liek$ consists of diagonal and
$\Liep$ consists of off-diagonal matrices. Then we consider the
$\Liesu(2)$-valued differential form
\BEQ
\talpha=\tF\inv\diff\tF.
\EEQ
It decomposes relative to the Cartan decomposition
$\talpha=\talpha_k+\talpha_p$. The $\Liep$-valued one-form $\talpha_p$
can be further decomposed into a holomorphic and an antiholomorphic
part $\talpha_p=\talpha_p^\prime\diffz+\talpha_p^\dprime\diffzbar$,
which gives
\BEQ
\talpha=\tF\inv\diff\tF=\talpha_p^\prime\diffz+\talpha_k
+\talpha^\dprime_p\diffzbar.
\EEQ
For each $\lambda\in S^1$ we define the following $\Liesu(2)$-valued one-form:
\BEQ \label{alphadef}
\alpha=\lambda\inv\talpha_p^\prime\diffz+\talpha_k
+\lambda\talpha_p^\dprime\diffzbar.
\EEQ
Because of~\cite{DoPeWu:1} and \cite{RuVi:1} we know, that $\tF$ is
the frame of a CMC-immersion iff $\alpha$ is integrable. In this case
we can define the extended frame
$F:\FD\times S^1\longrightarrow\LieSU(2)$ by 
\BEQ \label{Finitial}
\alpha=F\inv\diff F,\kern3cm F(0,\lambda)=I,\;\lambda\in S^1.
\EEQ
Then 
\BEQ \label{initial}
F(z,\zquer,1)=\tF(z,\zquer)
\EEQ
by the uniqueness of the initial value problem.
We will interpret $F$ in terms of loop groups in the next section.

\newsection \label{loopinterpretation}
Let $\Lambda\LieG_\sigma$, $\LieG=\LieSL(2,\Bcc)$ or $\LieG=\LieSU(2)$, 
denote the group of 
smooth maps $g(\lambda)$ from $S^1$ to $\LieG$, which satisfy the twisting
condition
\BEQ \label{DPWtwistcond}
g(-\lambda)=\sigma(g(\lambda)),
\EEQ
where $\sigma:\LieSU(2)\rightarrow\LieSU(2)$ is the group automorphism
of order $2$, which is given by conjugation with the Pauli matrix 
\BEQ
\sigma_3=\tmatrix100{-1}.
\EEQ
The Lie algebras of these groups, which we denote by
$\Lambda\Lieg_\sigma$, where $\Lieg$ is the Lie algebra of $\LieG$,
consist of maps $x:S^1\rightarrow\Lieg$, which satisfy a similar
twisting condition as the group elements
\BEQ
x(-\lambda)=\sigma_3 x(\lambda)\sigma_3.
\EEQ
In order to make these loop groups complex Banach Lie groups, we equip them,
as in \cite{DoPeWu:1}, with some $H^s$-topology for $s>{1\over2}$.
Elements of these twisted loop
groups are matrices with off-diagonal entries which are odd functions, and
diagonal entries which are even functions in the parameter $\lambda$.
All entries are in the Banach algebra ${\cal A}$ of $H^s$-smooth functions.

Furthermore, we will use the following subgroups of
$\Lambda\LieSL(2,\Bcc)_\sigma$: 
Let $\LieB$ be a subgroup of $\LieSL(2,\Bcc)$ and 
$\Lambda^+_B\LieSL(2,\Bcc)_\sigma$ be the group of maps in
$\Lambda\LieSL(2,\Bcc)_\sigma$, which can be extended to holomorphic maps on
the open unit circle and take values in $\LieB$ at $\lambda=0$.
Analogously, let $\Lambda^-_B\LieSL(2,\Bcc)_\sigma$ 
be the group of maps in $\Lambda\LieSL(2,\Bcc)_\sigma$, which can be extended
to the outside of the unit circle in $\CPE=\Bcc\cup\{\infty\}$ and
take values in $\LieB$ at $\lambda=\infty$. If $\LieB=\{I\}$ (based loops) 
we write the subscript $\ast$ instead of $\LieB$, if
$\LieB=\LieSL(2,\Bcc)$ we omit the subscript for $\Lambda$ entirely.

Corresponding to these subgroups we define Lie subalgebras of
$\Lambda\Liesl(2,\Bcc)_\sigma$: We will only consider
$\Lambda^+\Liesl(2,\Bcc)_\sigma$, which consists of maps $x(\lambda)$
that can be continued holomorphically to the unit circle, and
$\Lambda^-_\ast\Liesl(2,\Bcc)_\sigma$, which consists of maps, which
can be continued holomorphically to the outside of the unit circle and
become the zero matrix at $\lambda=\infty$.

We quote the following result from~\cite{DoPeWu:1}:

\separate (i) For each solvable subgroup $\LieB$ of $\LieSL(2,\Bcc)$
satisfying $\LieSU(2)\cdot\LieB=\LieSL(2,\Bcc)$ and
$\LieSU(2)\cap\LieB=\{I\}$, multiplication 
$$
\Lambda\LieSU(2)_\sigma\times\Lambda^+_B\LieSL(2,\Bcc)_\sigma
\longrightarrow\Lambda\LieSL(2,\Bcc)_\sigma
$$
is a diffeomorphism onto.
The associated splitting
\BEQ \label{Iwasawa}
g=F g_+
\EEQ
of an element $g$ of $\Lambda\LieSL(2,\Bcc)_\sigma$, s.t.\
$F\in\Lambda\LieSU(2)_\sigma$ and $g_+\in\Lambda^+_B\LieSL(2,\Bcc)_\sigma$
will be called Iwasawa decomposition (relative to $\LieB$). If not
specified otherwise, we will always choose $\LieB$ to be the set of upper
triangular matrices in $\LieSL(2,\Bcc)$ which have real, positive
diagonal entries.

\separate (ii) Multiplication 
$$
\Lambda^-_\ast\LieSL(2,\Bcc)_\sigma\times\Lambda^+\LieSL(2,\Bcc)_\sigma
\longrightarrow\Lambda\LieSL(2,\Bcc)_\sigma
$$
is a diffeomorphism onto the open and dense subset 
$\Lambda^-_\ast\LieSL(2,\Bcc)_\sigma\cdot\Lambda^+\LieSL(2,\Bcc)_\sigma$
of $\Lambda\LieSL(2,\Bcc)_\sigma$, called the ``big cell'' \cite{SeWi:1}.
The associated splitting
\BEQ
g=g_-g_+
\EEQ
of an element $g$ of the big cell, where
$g_-\in\Lambda^-_\ast\LieSL(2,\Bcc)_\sigma$ and
$g_+\in\Lambda^+\LieSL(2,\Bcc)_\sigma$, will be called Birkhoff factorization.

Using the definitions above, it follows from~\bref{alphadef}, that
$\alpha$ is a one-form taking values in $\Lambda\Liesu(2)_\sigma$. By
the definition~\bref{Finitial} of the extended frame, we see, that $F$
can be interpreted as a map from $\FD$ into $\Lambda\LieSU(2)_\sigma$.

\mytheorem{} {\em
The extended frame $F$ transforms under $g\in\Aut_\Psi\FD$ as follows:
\BEQ \label{exttransform}
(F\circ g)(z,\zquer,\lambda)=\chi(g,\lambda)F(z,\zquer,\lambda)k(g,z,\zquer),
\EEQ
where $\chi(g,\lambda)\in\Lambda\LieSU(2)_\sigma$.
In addition, there exists a map
\BEQ
\epsilon:\Aut_\Psi\FD\times\Aut_\Psi\FD\rightarrow\{I,-I\}\subset\LieSU(2),
\EEQ
s.t.\ for each $\lambda\in S^1$, the map
$\chi(\cdot,\lambda):\Aut_\Psi\FD\rightarrow\LieSU(2)$ is a projective group
homomorphism with cocycle $\epsilon$:
\BEQ \label{chihomo}
\chi(g_2\circ g_1,\lambda)=\epsilon(g_2,g_1)
\chi(g_2,\lambda)\chi(g_1,\lambda),
\EEQ
and $k(\cdot,z,\zquer)$ satisfies the projective cocycle relation
\BEQ \label{kcocyclecond}
k(g_2\circ g_1,z,\zquer)
=\epsilon(g_2,g_1)k(g_2,g_1(z),\overline{g_1(z)})k(g_1,z,\zquer).
\EEQ
}

\Proof
The differential $\talpha=\tF\inv\diff\tF$
transforms under an automorphism $g\in\Aut_\Psi\FD$ as
\BEQ
g_\ast\talpha=k\inv\talpha k+k\inv\diff k.
\EEQ
Using the Cartan decomposition above we get 
\BEA \label{Cartandecomp}
g_\ast\talpha_k & = & k\inv\talpha_k k+k\inv\diff k,\nonumber\\
g_\ast\talpha_p & = & k\inv\talpha_p k.
\EEA
By further decomposing $\talpha_p$ into a holomorphic and an
antiholomorphic part, we get
\BEA \label{holantiholdecomp}
(\talpha\circ g)_p^\prime \partial_zg & = & k\inv\talpha_p^\prime k,\nonumber\\
(\talpha\circ g)_p^\dprime\overline{\partial_zg} 
& = & k\inv\talpha_p^\dprime k.
\EEA
If we define $\cF(z,\zquer,\lambda)=F(z,\zquer,\lambda)k(g,z,\zquer)$, then
\BEQ
\calpha=\cF\inv\diff\cF=\lambda\inv
k\inv\talpha_p^\prime k\diffz+k\inv\talpha_k k+\lambda
k\inv\talpha_p^\dprime k\diffzbar+k\inv\diff k.
\EEQ
Therefore, by Eqs.~\bref{Cartandecomp}, \bref{holantiholdecomp} and
the definition~\bref{alphadef},
\BEQ
g_\ast\alpha=\calpha
\EEQ
and $\cF$ and $F\circ g$ are equal up to left
multiplication by a unitary $z$-independent matrix $\chi(g,\lambda)$,
\BEQ
(F\circ g)(z,\zquer,\lambda)=\chi(g,\lambda)F(z,\zquer,\lambda)k(g,z,\zquer).
\EEQ
The matrix $\chi(g,\lambda)$ is fixed by the initial 
condition~\bref{initial} as
\BEQ \label{chiinitial}
\chi(g,\lambda)=F(g(0),\lambda)k(g,0)\inv\in\Lambda\LieSU(2)_\sigma.
\EEQ
This shows, that $\chi(g,\lambda)k(g,0)$ is uniquely determined. Since
$k$ satisfies~\bref{kcocycle},
we can find for each pair of automorphisms $g_1$
and $g_2$ in $\Aut_\Psi\FD$ an $\epsilon(g_2,g_1)\in\{I,-I\}$, s.t.\
Eq.~\bref{kcocyclecond} holds. This in turn gives Eq.~\bref{chihomo}, since
\BEA
\chi(g_2\circ g_1,\lambda)k(g_2\circ g_1,0) & = &
F((g_2\circ g_1)(0)) \nonumber\\
& = & \chi(g_2,\lambda)\chi(g_1,\lambda)k(g_1,0)k(g_2,g_1(0)) \nonumber\\
& = & \epsilon(g_2,g_1)\chi(g_2,\lambda)\chi(g_1,\lambda)k(g_2\circ g_1,0).
\EEA
\mbox{\kern1cm}\QED

\newsection \label{Symsection}
Let us now look at Sym's formula
\BEQ \label{Sym}
J(\Psi_\lambda)=\pder{}{\theta}F\cdot F\inv
+{i\over2}F\sigma_3 F\inv,\;\lambda=e^{i\theta},
\EEQ
which gives for each extended frame $F$ an associated family (see
Section~\ref{associatedfamily}) of CMC-immersions 
$\Psi_\lambda:\FD\rightarrow\threespace$ in the spinor representation.
Remember, that we have chosen $H=-\frac{1}{2}$.

The initial condition in~\bref{Finitial} implies (see~\cite[Remark
A.7]{DoHa:1}) an initial condition for $\Psi_\lambda$:
\BEQ \label{Phiinitial}
\Psi_\lambda(z=0)=-e_3,
\EEQ
where $e_3$ is the unit vector mapped to $-{i\over2}\sigma_3$ by the
spinor representation. We will in the following always assume, that
the members of the associated family are normalized by~\bref{Phiinitial}.
Note, that then Eq.~\bref{assgroup} holds for the normalized surfaces.

Using Eq.~\bref{exttransform}, it is easy to see that the 
family of immersions $\Psi_\lambda$ transforms under $g\in\Aut_\Psi\FD$ as 
\BEA
J((\Psi_\lambda\circ g)(z)) & = & 
\pder{}{\theta}\chi\cdot\chi\inv
+\chi \pder{}{\theta}F\cdot F\inv \chi\inv
+{i\over2}\chi F\sigma_3 F\inv\chi\inv\nonumber\\
& = & \chi(J(\Psi_\lambda(z)))\chi\inv
+\pder{}{\theta}\chi\cdot\chi\inv.
\label{Symtrafo}
\EEA
The second term in Eq.~\bref{Symtrafo} is a translation in $\threespace$,
the first term describes a rotation of the initial surface. For
$\lambda=1$ we have $\Psi_1=\Psi$, $\Psi$ being the original
immersion. Therefore $R_{\tT}$, $\tT=(R_{\tT},t_{\tT})=\psi(g)$, 
acts on vectors in the spinor representation by 
conjugation with $\chi(g,\lambda=1)$ and the translation $t_T$ acts by 
adding 
$\pder{}{\theta}|_{\theta=0}\chi\cdot\chi\inv(g,\lambda=1)$.

\myprop{} {\em $\Aut_{\Psi_\lambda}\FD=\Aut_\Psi\FD$ for all
$\lambda\in S^1$.
}

\Proof From Eq.~\bref{Symtrafo} we see that for 
$g\in\Aut_\Psi\FD$ and arbitrary $\lambda\in S^1$ we have
\BEQ
\Psi_\lambda\circ g=\tT_\lambda\circ\Psi
\EEQ
where $\tT_\lambda\in\OAff(\threespace)$ is given in the spinor
representation by conjugation with $\chi(g,\lambda)$ and subsequent
addition of
$\pder{}{\theta}\chi\cdot\chi\inv(g,\lambda)$.
This implies $\Aut_\Psi\FD\subset\Aut_{\Psi_\lambda}\FD$. Since, by
Eq.~\bref{assgroup}, the surface $\Psi=\Psi_1$ is in the 
associated family of $\Psi_{\lambda}$ for each $\lambda\in S^1$,
we get also
\BEQ
\Aut_{\Psi_\lambda}\FD
=\Aut_{(\Psi_\lambda)_{\lambda\inv}}\FD\subset\Aut_\Psi\FD,
\EEQ
which finishes the proof.
\QED

\separate Using Proposition~\ref{Symsection} and
Lemma~\ref{SmythDel2}, we get the following

\mycorollary{} {\em 
Let $(M,\Phi)$ be a complete, nonspherical CMC-surface with simply
connected cover $(\FD,\Psi)$.
If $\Aut\Psi(\FD)$ contains a one-parameter group of proper 
Euclidean motions, then $(\FD,\Psi)$ is in the associated family of a
Delaunay surface.
}

\Proof 
If $(\FD,\Psi)$ is in the associated family of a Delaunay surface,
then by Lemma~\ref{SmythDel2} and~Proposition~\ref{Symsection},
$\Aut_\Psi\FD$ and therefore, by Theorem~\ref{mapdefs}, 
also $\Aut\Psi(\FD)$ contains a one-parameter group.

Conversely, if $\Aut\Psi(\FD)$ contains a one-parameter group, then, by 
Proposition~\ref{nondiscrete} and Corollary~\ref{secfund}, 
the group $\Aut_\Psi\FD$ also contains a
one-parameter group of self-isometries of $(\FD,\Psi)$.
By Theorem~\ref{DelundSmyth}, only those
surfaces which are isometric to Smyth or
Delaunay surfaces can have a one parameter group of self-isometries.
{}From Lemma~\ref{SmythDel2}
we know, that for Smyth surfaces $\Aut\Psi(\FD)$ is discrete.
This shows, that $(\FD,\Psi)$ is isometric to a Delaunay surface,
which, by Lemma~\ref{associatedfamily}, proves the claim.
\QED

\separate 
For each CMC-immersion $\Psi_\lambda$ in the associated family of
$\Psi$ we get, as in Section~\ref{mapdefs}, an analytic group homomorphism
\BEQ
\psi_\lambda:\Aut_{\Psi_\lambda}\FD=\Aut_\Psi\FD
\rightarrow\Aut\Psi_\lambda(\FD).
\EEQ
Since by Corollary~\ref{associatedfamily}, 
with $(\FD,\Psi)$ also $(\FD,\Psi_\lambda)$ is complete, we get,
by Corollary~\ref{immersedautom}, that $\psi_\lambda$ is surjective.
If we write the canonical covering $P:\LieSU(2)\rightarrow\LieSO(3)$
using the spinor representation $J$ of $\Brr^3$,
\BEQ
P(A)(x)=\Ad(A)(J(x)),\kern1cm A\in\LieSU(2),\kern0.5cm x\in\threespace,
\EEQ
we get from Eq.~\bref{Symtrafo} for fixed $\lambda\in S^1$:
\BEQ
\psi_\lambda(g)=\tT_\lambda=(R_{\tT_\lambda},t_{\tT_\lambda}),
\EEQ
where for $x\in\threespace$,
\BEQ \label{rottrans}
J(R_{\tT_\lambda}(x))=\Ad(\chi(g,\lambda))(J(x)),\kern3cm 
J(t_{\tT_\lambda}(x))=J(x)+\pder{}{\theta}\chi\cdot\chi\inv(g,\lambda).
\EEQ
Since, by~Theorem~\ref{mapdefs}, $\psi_\lambda$ is analytic, we get

\mylemma{} {\em
For each $\lambda\in S^1$,
\BEQ
\Ad(\chi(\cdot,\lambda)):\Aut_\Psi\FD\longrightarrow\Aut(\Liesu(2))
\EEQ
is a continuous group homomorphism, and the map
\BEQ
\pder{}{\theta}\chi\cdot\chi\inv(\cdot,\lambda):\Aut_\Psi\FD\longrightarrow
\Liesu(2),\;\lambda=e^{i\theta},
\EEQ
is continuous.
}

\Proof 
Since the decomposition~\bref{congdecomp} of a proper Euclidean motion
is unique, we get, that with $\tT_\lambda=\psi_\lambda(g)$ 
also the rotational and
translational parts, $R_{\tT_\lambda}$ and $t_{\tT_\lambda}$, depend
continuously on $g\in\Aut_{\Psi_\lambda}\FD=\Aut_\Psi\FD$. 
Using Eq.~\bref{rottrans}, we get the lemma.
\QED

\myremark{}
Let us define a Riemann surface $M_\lambda=\Ker\psi_\lambda\setminus\FD$
as in Proposition~\ref{groups}. Let $\pi_\lambda:\FD\rightarrow M_\lambda$ 
denote the natural covering map and let $\bpi_\lambda$ be defined as
in Definition~\ref{mapdefs}. Then we have
\BEQ
\Ker\bpi_\lambda=\Ker\psi_\lambda
\EEQ
for all $\lambda\in S^1$. This shows, that if we define $M_\lambda$ as
above, then Condition~\bref{simplycovered} is satisfied for all
members of the associated family. However, it should be noted, that
$\Ker\psi_\lambda$ and therefore also $M_\lambda$ can vary in a
very complicated way with $\lambda$.

\newsection \label{meropotentials}
Now we apply the Birkhoff splitting $F=g_-g_+$ for
$z\in\FD\setminus\FS$, where $\FS$ is the discrete set of points in
$\FD$, where $F$ cannot be split. Then $g_-$ can be continued
meromorphically to $z\in\FD$. For a detailed account see
\cite[Section~4]{DoPeWu:1}, especially the remarks following Theorem~4.10.

By Eq.~\bref{exttransform} we get for $g\in\Aut_\Psi\FD$
\BEQ \label{gmintrafo}
(g_-\circ g)(z,\lambda)=\chi(g,\lambda)g_-(z,\lambda)
p_+(g,z,\lambda),
\EEQ
where
\BEQ
p_+=g_+(z,\zquer),\lambda)k(g,z,\zquer)
((g_+\circ g)(z,\zquer,\lambda))\inv
\EEQ
is determined by $g$ up to a sign.
Since $k$ satisfies Eq.~\bref{kcocyclecond} we have for
$p_+(g,z,\lambda)$:
\BEQ
p_+(g_2\circ g_1,z,\lambda)=\epsilon(g_2,g_1) 
p_+(g_2,g_1(z),\lambda)p_+(g_1,z,\lambda),
\EEQ
where
$\epsilon:\Aut_\Psi\FD\times\Aut_\Psi\FD\rightarrow\{I,-I\}\subset\LieSU(2)$
is defined as in Theorem~\ref{loopinterpretation}.

The meromorphic potential $\xi$ is defined by $\xi=g_-\inv\diff g_-$.
It is a $\Lambda^-\Liesl(2,\Bcc)_\sigma$-valued one-form with only a
$\lambda\inv$-coefficient. Therefore, it is of the form
\BEQ
\xi=\lambda\inv\tmatrix0f{E\over f}0\diffz,
\EEQ
where $f$ is a meromorphic function and $E\diffz^2$ is the holomorphic
Hopf differential.
{}From Eq.~\bref{gmintrafo} we obtain, that $\xi$ 
transforms under the automorphism $g\in\Aut_\Psi\FD$ as
\BEQ \label{xitransform}
g_\ast\xi(z,\lambda)=(\xi\circ g)(z,\lambda)g^\prime(z)
=p_+\inv\xi p_++p_+\inv\diffp_+.
\EEQ
Elements of $\Aut_\Psi\FD$ therefore act as 
gauge transformations on the meromorphic potential.

\myremark{}
By the normalization of $F$ we have
\BEQ
(F\circ g)(0,\lambda)=\chi(g,\lambda)k(g,0).
\EEQ
If $g(0)$ is not in the singular set $\FS$, then, since
$g_-(0,\lambda)=I$, we can evaluate~\bref{gmintrafo} at $z=0$. This
shows $g_-(g(0),\lambda)=\chi(g,\lambda)p_+(g,0,\lambda)$. In
particular $\chi(g,\lambda)$ is splittable. Therefore,
\BEQ
g_-(g(0),\lambda)=\chi_-(g,\lambda),
\EEQ
where
\BEQ
\chi=\chi_-\chi_+
\EEQ
is the Birkhoff splitting of $\chi$.
{}From Eq.~\bref{gmintrafo} it follows
\BEQ \label{gdressing}
\chi_-\inv(g,\lambda)(g_-\circ g)(z)
=\chi_+(g,\lambda)g_-(z,\lambda)p_+(g,z,\lambda).
\EEQ
Now
\BEQ
g^\chi_-(z,\lambda)=\chi_-\inv(g,\lambda)(g_-\circ g)(z,\lambda)
\EEQ
is normalized at $z=0$. Thus Eq.~\bref{gdressing} can be looked at
as a dressing transformation from $g_-$ to $g^\chi_-$, a new surface.

\newsection \label{automorphism}
The procedure presented in~\cite{DoPeWu:1} and used in this section
characterizes CMC-surfaces in terms of the meromorphic potential $\xi$
(for smoothness questions see also~\cite{DoHa:1}). It is
therefore natural to attempt to characterize $g\in\Aut_\Psi\FD$ by
using equations only involving $\xi$. The obvious equation to start
with is~\bref{xitransform}. Thus, let us assume that to some
$g\in\Aut\FD$ we can find some
$p_+=p_+(g,z,\lambda)\in\Lambda^+\LieSL(2,\Bcc)_\sigma$, s.t.\
\BEQ \label{vierzehn1}
(g_\ast\xi)(z,\lambda)=p_+\inv\xi p_++p_+\inv\diff p_+.
\EEQ
Note, that every such $p_+$ is automatically meromorphic in $z$. Next
we solve
\BEQ \label{vierzehn2}
g_-\inv\diff g_-=\xi,\kern3cm g_-(0,\lambda)=I
\EEQ
and consider $h=g_-p_+$. In view of Eq.~\bref{vierzehn1}
and~\bref{vierzehn2} it is straightforward to verify that $h\inv\diff
h=g_\ast\xi$ holds. On the other hand, $\hh=g_-\circ g$ also satisfies
$\hh\inv\diff\hh=g_\ast\xi$. Therefore, 
\BEQ \label{vierzehn3}
(g_-\circ
g)(z,\lambda)=\sigma(g,\lambda)g_-(z,\lambda)p_+(g,z,\lambda).
\EEQ
Note, that here $\sigma\in\Lambda\LieSL(2,\Bcc)_\sigma$ does not depend on $z$.

In order to be able to conclude, that $g\in\Aut_\Psi\FD$, we need to
obtain Eq.~\bref{exttransform}, since the discussion in
Section~\ref{Symsection} shows that Eq.~\bref{exttransform} implies
$g\in\Aut_\Psi\FD$.
The relation between~\bref{exttransform} and~\bref{gmintrafo} is
obtained via the Birkhoff splitting
\BEQ
(g_-\circ g)\cdot(g_+\circ g)=\chi g_-g_+k.
\EEQ
We therefore need 
\BEQ
\chi\inv\sigma g_-=g_-q_+,\kern1cm
q_+=q_+(g,z,\lambda)\in\Lambda^+\LieSL(2,\Bcc)_\sigma.
\EEQ
But then, setting $z=0$, we obtain
\BEQ
\chi\inv\sigma\in\Lambda^+\LieSL(2,\Bcc)_\sigma.
\EEQ
Setting $h_+=\chi\inv\sigma$ we thus obtain $h_+g_-=g_-q_+$. This
shows
\BEQ
\sigma g_-=\chi h_+g_-=\chi g_-q_+,
\EEQ
whence
\BEQ
g_-\circ g=\sigma g_-p_+=\chi g_-q_+p_+.
\EEQ
Altogether, we have shown

\myprop{} {\em
Let $g\in\Aut\FD$. Then $g\in\Aut_\Psi\FD$ iff there exists some
$p_+=p_+(g,z,\lambda)\in\Lambda^+\LieSL(2,\Bcc)_\sigma$ s.t.\
\BEQ \label{ppeq}
g_\ast\xi=p_+\inv\xi p_++p_+\inv\diff p_+
\EEQ
and the $z$-independent matrix function $\sigma$ associated with $p_+$
by~\bref{vierzehn3} is unitary.
}

\myremark{}
We have seen in the proof that, starting from some $p_+$ satisfying
Eq.~\bref{vierzehn1}, we may have to modify the corresponding $\sigma$
by some $h_+$ s.t.\ $\chi=\sigma h_+\inv$ becomes unitary. In this
case $h_+$ satisfies $h_+g_-=g_-q_+$, i.e., $h_+$ is in the isotropy
group of the dressing action on $g_-$. 

\newsection \label{cyldressed}
To illustrate the discussion of the last section we look
at the dressing orbit of the standard cylinder.

The meromorphic potential of the cylinder can be chosen to be
\BEQ
\xi_0=\lambda\inv A\diffz,
\EEQ
where
\BEQ
A=\tmatrix0110
\EEQ
is constant.
We can integrate $\xi_0$ without problems:
\BEQ
g^0_-=\exp\left(\lambda\inv zA\right),
\EEQ
and, since $A$ commutes with its adjoint matrix $\bA=A$, also the Iwasawa
splitting doesn't pose any problems:
\BEQ \label{standardcylframe}
F=\exp\left((\lambda\inv z-\lambda \zquer)A\right).
\EEQ
We will call the cylinder, which is generated by applying
Sym's formula to Eq.~\bref{standardcylframe}, the {\em standard cylinder}.
For surfaces in the dressing orbit of the standard cylinder we therefore have
\BEQ \label{gmindressdef}
g_-=(h_+g^0_-)_-=h_+g^0_-g_+\inv,\kern3cm\xi=g_-\inv\diffg_-,
\EEQ
with $h_+\in\Lambda^+\LieSL(2,\Bcc)_\sigma$.
For these $g_-$ 
each translation $z\mapsto z+q$, $q\in\Bcc$, in $\FD$ can be
implemented by a transformation of type~\bref{gmintrafo},
\BEQ \label{dresstrans}
g_-(z+q,\lambda)=Q\cdot g_-(z,\lambda) r_+(q,z,\lambda),
\EEQ
with the $z$-independent matrix
\BEQ \label{Qdef}
Q=h_+e^{\lambda\inv qA}h_+\inv
\EEQ
and 
\BEQ
r_+(q,z,\lambda)=g_+(z,\lambda)g_+\inv(z+q,\lambda).
\EEQ
To illustrate Proposition~\ref{automorphism}, we consider the standard
cylinder itself,
i.e., $h_+=I$ and $g_-=e^{\lambda\inv zA}$. Then for every $q\in\Bcc$,
$Q=e^{\lambda\inv qA}$. Moreover, $Qe^{-\lambda\qquer A}$ is unitary and
$e^{-\lambda\qquer A}g_-=g_-p_+$. 
Therefore, Proposition~\ref{automorphism} shows
that every translation $z\mapsto z+q$, $q\in\Bcc$, is in
$\Aut_\Psi\FD$. But for general $h_+$ the matrix 
$Q$ cannot be modified in the sense of
Section~\ref{automorphism} to yield a unitary $\chi$.

In fact we can prove the following

\mytheorem{} {\em Other than cylinders there are no surfaces with
translational symmetry in the dressing orbit of the standard cylinder.
}

\underline{Proof:} Let $(M,\Phi)$ be a CMC-surface in the dressing
orbit of the standard cylinder with simply connected cover $(\FD,\Psi)$. 
Let $F:\FD\rightarrow\Lambda\LieSU(2)_\sigma$ be its extended frame
and let $g_-:\FD\rightarrow\Lambda^-\LieSL(2,\Bcc)_\sigma$ 
be defined by the Birkhoff splitting $F=g_-g_+$.
We assume, that $\Aut_\Psi\FD$ contains a translation $z\mapsto z+q$,
$q\in\cstar=\Bcc\setminus\{0\}$.

By Proposition~\ref{automorphism}, a translation $z\mapsto z+q$, 
$q\in\cstar$, is
in $\Aut_\Psi\FD$ iff $g_-$ satisfies
\BEQ
g_-(z+q,\lambda)=\chi(q,\lambda)g_-(z,\lambda)p_+(q,z,\lambda)
\EEQ
where $\chi(q,\lambda)\in\Lambda\LieSU(2)_\sigma$. 
By Equations~\bref{dresstrans} and \bref{Qdef}, this is equivalent to
\BEQ \label{periodeq}
Q(q,\lambda) g_-(z,\lambda)r_+(q,z,\lambda)=
\chi(q,\lambda)g_-(z,\lambda)p_+(q,z,\lambda),
\EEQ
where $Q=h_+e^{\lambda\inv qA}h_+\inv$. 
Setting $z=0$ in Eq.~\bref{periodeq}, it follows that $Q=\chi R_+$ with
$R_+=p_+(q,0,\lambda)r_+(q,0,\lambda)\inv\in\Lambda^+\LieSL(2,\Bcc)$.
Therefore, $g_-$ is invariant under
the dressing with $R_+$. This implies with Eq.~\bref{gmindressdef}
that for
\BEQ
w_+(q,\lambda,z)=g_+\inv p_+r_+\inv g_+\in\Lambda^+\LieSL(2,\Bcc)_\sigma,
\EEQ
we get
\BEQ \label{Rpeq}
R_+=g_-p_+r_+\inv g_-\inv=h_+e^{\lambda\inv zA}w_+ e^{-\lambda\inv zA}h_+\inv,
\EEQ
or equivalently
\BEQ
\chi=h_+e^{\lambda\inv (q+z)A}w_+\inv e^{-\lambda\inv zA}h_+\inv.
\EEQ
{}From Eq.~\bref{Rpeq} we get the condition
\BEQ \label{wpcond}
e^{\lambda\inv zA}w_+e^{-\lambda\inv zA}\in\Lambda^+\LieSL(2,\Bcc)_\sigma
\EEQ
for all $z\in\FD$. 
If we conjugate the expression in Eq.~\bref{wpcond} with 
\BEQ
D=\frac{1}{\sqrt{2}}\tmatrix11{-1}1,
\EEQ
then the exponent $DAD\inv=\sigma_3$ becomes diagonal and
with
\BEQ
\tw=DwD\inv=\tmatrix{\tw_a}{\tw_b}{\tw_c}{\tw_d},
\EEQ
Condition~\bref{wpcond} is equivalent to
\BEQ \label{blabla}
e^{\lambda\inv z\sigma_3}\tw e^{-\lambda\inv z\sigma_3}=
\tmatrix{\tw_a}{e^{2\lambda\inv z}\tw_b}{e^{-2\lambda\inv z}\tw_c}{\tw_d}
\EEQ
having a holomorphic extension to the unit disk for all
$z\in\FD\setminus\FS$. Since for every fixed $z\in\FD\setminus\FS$ the
off-diagonal entries of the r.h.s.\ of~\bref{blabla} vanish
identically or have an essential singularity at $\lambda=0$, this
implies that $\tw_b=\tw_c=0$. This
shows that $\tw$ commutes with $\sigma_3$. Thus $w$ commutes with $A$. 
We therefore have
\BEQ
\chi=h_+ H h_+\inv
\EEQ
with
\BEQ \label{Habform}
H(\lambda)=e^{\lambda\inv qA}w_+\inv(\lambda)=
\alpha(\lambda) I+\beta(\lambda) A,
\kern1cm\alpha(\lambda),\beta(\lambda)\in\Bcc.
\EEQ
Since $\det H=1$ we get also
\BEQ \label{alpha2beta2}
\alpha^2-\beta^2=1.
\EEQ
It follows from the unitarity of $\chi(q,\lambda)$ that
\BEQ \label{HPrel}
\Hqt=PH\inv P\inv,\kern2cm P=\hpqt h_+.
\EEQ
If we write Eq.~\bref{HPrel} using Eq.~\bref{Habform}, this gives
\BEQ \label{PAPEq}
(\alpha-\alphaquer)I=\beta PAP\inv+\betaquer A.
\EEQ
By taking the trace we get
\BEQ \label{alphareal}
\alpha=\alphaquer,
\EEQ
i.e., $\alpha$ is real, and
\BEQ \label{blabla2}
\betaquer AP=-\beta PA.
\EEQ
Using~\bref{alpha2beta2}, we get from~\bref{alphareal}
\BEQ \label{beta2real}
\beta^2=\betaquer^2.
\EEQ
{}From this and~\bref{blabla2} it follows,
that either 1) $\beta=\betaquer\equiv0$
or 2) $\beta=-\betaquer\not\equiv0$ and $[A,P]=0$ or 3)
$\beta=\betaquer\not\equiv0$ and $AP+PA=0$. 

In the first case, we get by~\bref{alpha2beta2}, that
$\alpha=\pm1$, i.e., $H=\chi=\pm I$. But this implies with
Eq.~\bref{Habform}, that $q=0$, which contradicts the assumption that
$q\in\cstar$.

In the second case, it follows from $[A,P]=0$, that $S=h_+Ah_+\inv$ is unitary.
The matrix $h_+Ah_+\inv$ can then be continued holomorphically to
$\CPE$, which is only possible if 
\BEQ \label{hpAcommute}
h_+Ah_+\inv=h_0Ah_0\inv,
\EEQ
where $h_0$ is a $\lambda$-independent diagonal matrix.
For the diagonal matrix $h_0$ the
unitarity of $h_0Ah_0\inv$ amounts to the unimodularity of the
diagonal entries. Therefore, the matrix $h_0$ is of the form
\BEQ
h_0=\tmatrix{e^{i\phi}}00{e^{-i\phi}},\kern3cm\phi\in\Brr.
\EEQ
By Eq.~\bref{hpAcommute} dressing of the standard cylinder with $h_+$ is
therefore the same as dressing with a unitary matrix $h_0$.
By writing the dressing action as an action on the extended
frames, we see, that dressing with a unitary matrix amounts to a 
rotation in $\threespace$. Therefore, in the second case the
surface $(M,\Phi)$ is a cylinder.

In the third case, it follows from $AP+PA=0$, that
$S=h_+Ah_+\inv$ satisfies $S\Squer^\top=-S^2=-I$.
This is impossible since $S\Squer^\top$ is positive
definite. Therefore, the third case cannot occur.
\QED

\mycorollary{ 1} {\em 
Let $(M,\Phi)$ be obtained from the standard cylinder $(M_0,\Phi_0)$ by
dressing with $h_+\in\Lambda^+\LieSL(2,\Bcc)_\sigma$. If $(M,\Phi)$
admits a translational symmetry, then $(M,\Phi)$ is a cylinder and
$h_+=h_0H_+$, where $h_0$ is unitary and $\lambda$-independent,
and $H_+$ commutes with $A$.
}

\separate
Theorem~\ref{cyldressed} shows, that there are neither Delaunay
surfaces nor CMC-tori in the dressing orbit of the standard cylinder.
It therefore also proves:

\mycorollary{ 2} {\em
The dressing action does not act transitively on the set of
all CMC-surfaces without umbilics.
}

\separate
The method developed in Section~\ref{automorphism} can also be applied
to the more general dressing action of~\cite{DoWu:1}. In this case the
orbit of the standard cylinder under this action contains all surfaces of
finite type, therefore all CMC-tori
(see~\cite{PiSt:1}). In~\cite{DoHa:3} we will use this method to
derive the classification of CMC-tori given in \cite{Bo:1}
in a more geometric way.

\newsection
Since we are interested in CMC-surfaces with symmetry groups, we would
like to understand what it means for $\xi$ and $g\in\Aut\FD$ to admit
some $p_+$ satisfying Eq.~\bref{ppeq}.
Expanding $p_+=p_0+\lambda p_1+\ldots$ into powers of $\lambda$, we
see, that~\bref{ppeq} implies
\BEQ
g_\ast\xi=p_0\inv\xi p_0.
\EEQ
Evaluating this with $\xi=\lambda\inv\tmatrix0f{\frac{E}{f}}0\diffz$ we
obtain for every $g\in\Aut_\Psi\FD$
\BEQ
(E\circ g)(z)(g^\prime(z))^2=E(z),
\EEQ
which we derived already in Lemma~\ref{secfund}, and in addition
\BEQ \label{ftrafo}
(f\circ g)(z)g^\prime(z)=r^{-2}(g,z)f(z),
\EEQ
where $p_0=\diag(r,r\inv)$.
The function $r^{-2}(g,z)$ is uniquely determined by
$g\in\Aut_\Psi\FD$. Using Eq.~\bref{ftrafo}, it is easy to verify that 
$r^{-2}$ satisfies the cocycle condition
\BEQ \label{dreineunvier}
r^{-2}(g_2\circ g_1,z)=r^{-2}(g_2,g_1(z))r^{-2}(g_1,z),
\EEQ
for $g_1,g_2\in\Aut_\Psi\FD$ and $z\in\FD$.

\myremark{} 
1. If $\FD=\Bcc$ and $G\subset\Aut_\Psi\FD$ is a group of (pure)
translations, then $g^\prime(z)=1$ for all $g\in G$. Hence, $E$ is
invariant under $G$. E.g., for tori, $E$ is doubly periodic and
holomorphic. As a consequence, $E$ is constant. This just reproduces the
well known fact, that CMC-tori do not have any umbilics.

2. If $\FD$ is the unit disk, then $\Aut\FD\cong\LieSU(1,1)$ consists
of transformations of the form $g(z)=\frac{az+b}{cz+d}$. Hence,
$g^\prime(z)=(cz+d)^{-2}$, since $\det\tmatrix abcd=1$.
Therefore, in this case $E$ is an automorphic form of degree $4$ for
$\Aut_\Psi\FD$.

3. For $f$ the situation is unfortunately less straightforward. A
condition on $E$ and $f$ can only be derived, if one is able to
express $r^{-2}$ in terms of $E$ and $f$ and
$g\in\Aut_\Psi\FD$. This will be discussed in Chapter~\ref{unitary}.

4. We would also like to note, that for any $g\in\Aut\FD$ we can use
Eq.~\bref{ftrafo} as a definition of the function $s_0=r^{-2}$.

\newsection \label{complete}
If $M$ is a Riemann surface with simply connected cover $\FD=\Bcc$,
then by Lemma~\ref{Riemannian} 
it is conformally equivalent to the complex plane, the cylinder
or a torus. Of course, if $M$ is a torus, then $M$ is complete. If $M$
is the plane or a cylinder, it would be interesting to know, what 
it means for the meromorphic
potential, that $M$ with the induced metric is complete.
For the case that the universal cover of $M$ is the unit disk, we can show

\mytheorem{} {\em
Let $(M,\Phi)$ be a CMC-immersion with universal cover $\FD$ the unit
disk. Let $\xi=\lambda\inv\tmatrix0f{\frac{E}{f}}0\diffz$ be the meromorphic
potential for $(M,\Phi)$. Then if $M$ is complete, $\xi$ cannot be
extended meromorphically to any open set containing a point of the
unit circle.
}

\Proof 
If $\xi$ could be extended meromorphically to an open set $U$ which
contains a point on the unit circle, then $U$ contains
an open subset $U_1$, which contains a point on the unit circle,
and on which $\xi$ is holomorphic.

If we take a curve inside $U_1$ which connects a point inside $\FD$ with a
point on the unit circle, then $\xi$, and therefore also the 
surface, can be extended along this line. The function $u(z)$ defined
in Eq.~\bref{udef} satisfies \cite[Sect.~A.8]{DoHa:1}
\BEQ \label{metricandf}
w^2(z,\zquer)f(z)=\frac{1}{4}e^{u(z,\zquer)\over2}
\EEQ
where $w(z,\zquer)$ is a map which is defined and nonzero on
$U_1\cup(\FD\setminus\FS)$, the domain of
definition of $\xi$. Therefore, the line has finite length w.r.t.\ the 
metric on $M$, which contradicts the completeness of $M$.
\QED

\section{Unitary gauge condition} \label{unitary}\message{[unitary]}

In this chapter we will use the results of Chapter~\ref{DPWautom}
to derive conditions on the meromorphic potential for the associated
CMC-immersion to be symmetric.

\newsection \label{aundb}
Let
\BEQ
\xi=\lambda\inv\tmatrix0f{E\over f}0\diffz,\kern2cm
\hxi=\lambda\inv\tmatrix0{\hf}{\hE\over\hf}0\diffz
\EEQ
be defined on $\FD$, where $\FD$ is the unit disk or the complex plane.
We define $g_-$ and $\hg_-$ by
\BEA \label{arbitraryinitial}
g_-^\prime & = & g_-\xi,\kern2cm g_-(0,\lambda)=g_0(\lambda),\nonumber\\
\hg_-^\prime & = & \hg_-\hxi,\kern2cm\hg_-(0,\lambda)=\hg_0(\lambda),
\EEA
where we have chosen arbitrary initial conditions
$g_0,\hg_0\in\Lambda^-\LieSL(2,\Bcc)_\sigma$.
In this section we want to investigate, under which conditions $g_-$
and $\hg_-$ give the same CMC-immersion up to a proper Euclidean
motion.

If we split the unitary matrix $\gquer g_-=\rquer r_+$ with
$r_+\in\Lambda^+\LieSL(2,\Bcc)_\sigma$, we have a $\LieU(1)$-freedom
to choose the $\lambda^0$-factor $r_0$ of $r_+$. We can use this freedom to
make $r_0$ real and positive.

Using this, we let now $D$ and $\hD$ be the real valued functions defined by 
\BEA
(\gquer g_-)_{+0} & = & \exp(D\sigma_3),\nonumber\\
(\hgquer\hg_-)_{+0} & = & \exp(\hD\sigma_3),
\EEA
where $g_{+0}$ denotes the real $\lambda^0$-coefficient of
the positive splitting factor of $g\in\Lambda\LieSL(2,\Bcc)_\sigma$.

With the definitions above we get the following

\mytheorem{} {\em
Let\/ $\xi$ and $\hxi$ be meromorphic differentials which give under
the DPW construction, using the initial conditions in~\bref{arbitraryinitial},
complete CMC-immersions $(\FD,\Psi)$ and $(\FD,\hPsi)$ with mean
curvature $H=-\frac{1}{2}$. 
Then the following statements are equivalent

1. The CMC-immersions $(\FD,\Psi)$ and $(\FD,\hPsi)$ 
differ only by a proper Euclidean motion, i.e., $\hPsi=\tT\circ\Psi$ for
some $\tT\in\OAff(\threespace)$.

2. There is a gauge transformation 
\BEQ \label{gaugeeq}
\hg_-=\chi g_-p_+,
\EEQ
with $\chi\in\Lambda\LieSU(2)_\sigma$ $z$-independent,
$p_+=p_+(z,\lambda)\in\Lambda^+\LieSL(2,\Bcc)_\sigma$.

3. $E=\hE$ and
\BEQ \label{secondcond}
\exp(D-\hD)={|\hf|\over|f|}.
\EEQ
}

\Proof 
Let $F$ and $\tF$ denote the extended frames of $(\FD,\Psi)$ and
$(\FD,\hPsi)$ defined by the Iwasawa decomposition of $g_-$ and
$\hg_-$, respectively.

1.$\Rightarrow$2.: 
If $(\FD,\Psi)$ and $(\FD,\hPsi)$ are related by a proper Euclidean
motion, then for their frames we get
\BEQ
\hF(z,\zquer,\lambda=1)=UF(z,\zquer,\lambda=1)k(z,\zquer),
\EEQ
where $U\in\LieSU(2)$ and $k:\FD\rightarrow\LieU(1)$ (see
Section~\ref{33}). By the same arguments as in Section~\ref{Chistetig}
and~\ref{loopinterpretation} we get for the extended frames
\BEQ \label{hFFinth}
\hF(z,\zquer,\lambda)=\chi(\lambda)F(z,\zquer,\lambda)k(z,\zquer),
\EEQ
with $\chi\in\Lambda\LieSU(2)_\sigma$. Here it should be noted, that
the initial condition in~\bref{Finitial} is replaced by
\BEQ
F(0,\lambda)=F_0(\lambda),\kern3cm \hF(0,\lambda)=\hF_0(\lambda),
\EEQ
where $F_0$ and $\hF_0$ are the unitary parts of $g_0$ and $\hg_0$,
respectively, w.r.t.\ the Iwasawa decomposition.
Using the Birkhoff splitting, we get Eq.~\bref{gaugeeq}
from~\bref{hFFinth}, with $p_+\in\Lambda^+\LieSL(2,\Bcc)_\sigma$.

2.$\Rightarrow$3.: Starting from $g_-$, $\hg_-$ satisfying
Eq.~\bref{gaugeeq} we consider the associated meromorphic potentials
\BEQ
\xi(\lambda)=\lambda\inv\tmatrix0f{E\over f}0\diffz\kern1cm\mbox{and}\kern1cm 
\hxi(\lambda)=\lambda\inv\tmatrix0\hf{\hE\over\hf}0\diffz.
\EEQ
Let $\frac{1}{2}e^u\diffz\diffzbar$ and
$\frac{1}{2}e^{\hu}\diffz\diffzbar$ be the metrics of the surfaces
defined by $\xi$ and $\hxi$. Then, by \cite[Sect.~A.8]{DoHa:1}, 
$f$ and $\hf$ are given by
\BEQ \label{metricfromf}
f=\frac{1}{4}w_0^{-2}e^{u\over2},\kern5cm
\hf=\frac{1}{4}\hw_0^{-2}e^{\hu\over2}.
\EEQ
Here $w_0$ and $\hw_0$ are defined by the Iwasawa decomposition
\BEQ \label{giwasplit}
\hg_-=\hF\hh_+,\hskip2cm g_-=Fh_+,
\EEQ
which is chosen, s.t.\ $u$ and $\hu$ in Eq.~\bref{metricfromf} are real
(see \cite{DoPeWu:1,DoHa:1}).
They are the upper left entries of the $\lambda^0$-coefficients of
$h_+$ and $\hh_+$, respectively.

Using the definition of $\exp D$ and
\BEQ
\gquer g_-=\hquer h_+,
\EEQ
we see that
\BEQ
\exp D=w_0^2e^{i\theta},\;\theta(z)\,\mbox{real},
\EEQ
and similarly 
\BEQ
\exp \hD=\hw_0^2e^{i\htheta},\;\htheta(z)\,\mbox{real}.
\EEQ
Taking the quotient, we get
\BEQ \label{genDhD}
\left|{\hf\over
f}\right|=\left|{w_0\over\hw_0}\right|^2e^{\hu-u\over2}=
\exp(D-\hD)e^{\hu-u\over2}.
\EEQ
This is true for any two CMC-immersions.
Next we invoke Eq.~\bref{gaugeeq}. Since $\chi$ is unitary, an Iwasawa
splitting gives $\hF=\chi Fk$ with some $\lambda$-independent, unitary
$k=k(z,\zquer)$. Using again that $\chi$ is unitary, this implies that
the induced metrics are the same, whence~\bref{secondcond}. Moreover,
from Eq.~\bref{gaugeeq} we get
\BEQ
\hxi=p_+\inv\xi p_++p_+\inv\diff p_+=p_0\inv\xi p_0,
\EEQ
where $p_0(z)=p_+(z,\lambda=0)$. From this we get
\BEQ
\hE=-\lambda^2\det(\hxi)=-\lambda^2\det(\xi)=E.
\EEQ

3.$\Rightarrow$1.:
If Eq.~\bref{secondcond} holds, then, by Eq.~\bref{genDhD}, $u=\hu$,
i.e., the CMC-immersions $(\FD,\Psi)$ and $(\FD,\hPsi)$ have the same metric
and the same mean curvature.
Since they also have the same Hopf differential, they
differ, by Corollary~\ref{DGL}, only by a proper Euclidean motion in
$\threespace$.
\QED

\myremark{} $D$ and $\hD$ can also be defined in the Grassmannian 
representation. There the Birkhoff splitting is replaced by the block
matrix splitting
\BEQ
\tmatrix ABCD=\tmatrix{a_-}0c{d_-} \tmatrix {a_+}b0{d_+},
\EEQ
(see \cite{SeWi:1,DoNeSz:1}). 
The upper diagonal block $a_+$ in the second factor can be chosen to
be an upper triangular matrix with a positive real number $s$ in the
lower right corner. This fixes the splitting factors completely and an
easy calculation shows that $s^2=\exp(-D)$ and analogously for $\hD$.

Let us define the associated families $\{\Psi_\lambda\}_{\lambda\in S^1}$,
$\{\hPsi_\lambda\}_{\lambda\in S^1}$ of the CMC-immersions $\Psi$ and
$\hPsi$ using extended frames and Sym's formula as in 
Section~\ref{Symsection}. In particular, we normalize the immersions
$\Psi_\lambda$ and $\hPsi_\lambda$ for all $\lambda$ by Eq.~\bref{Phiinitial}.
{}From the proof of Proposition~\ref{Symsection}, we know, 
that with $\Psi$ and $\hPsi$ also
$\Psi_\lambda$ and $\hPsi_\lambda$, $\lambda\in S^1$, are related by a unique
proper Euclidean motion $\tT_\lambda$. It is a natural
question, under what conditions on the meromorphic potentials $\xi$ and
$\hxi$, this Euclidean motion does not depend on $\lambda$. In
this case, the corresponding members of the associated families differ
all by the same proper Euclidean motion. This problem is treated in
the following

\mycorollary{} {\em Let $\xi$ and $\hxi$ be meromorphic potentials
associated with CMC-immersions. Let $g_-$ and $\hg_-$ be solutions of
\BEQ \label{initialx}
g_-^\prime=g_-\xi,\,g_-(0,\lambda)=I,\kern3cm\hg_-^\prime=\hg_-\hxi,\,
\hg_-(0,\lambda)=I.
\EEQ
Assume also $E=\hE$. Then the following are equivalent:

a) For the functions $f$ and $\hf$ we have
\BEQ \label{conditionb}
\exp(D-\hD)={|\hf|\over |f|}.
\EEQ

b) There exists a $\lambda$-independent, unitary matrix $\chi$, s.t.\
\BEQ \label{conclusion}
\hg_-=\chi g_-\chi\inv.
\EEQ

c) There exists a constant $c\in\Bcc$, $|c|=1$, s.t.\ $\hf=cf$.

d) There exists a $\lambda$-independent, unitary matrix
$\chi\in\LieU(1)$, s.t.\ for the
associated families $\Psi_\lambda$ and $\hPsi_\lambda$ defined by $\xi$ and
$\hxi$ we have in the spinor representation
\BEQ
J(\hPsi_\lambda)=\chi J(\Psi_\lambda)\chi\inv.
\EEQ
}

\Proof 
a)$\Rightarrow$b) 
{}From Theorem~\ref{aundb} and Eq.~\bref{conditionb} we get, that 
\BEQ \label{gmhut}
\hg_-(z,\lambda)=\chi(\lambda)g_-(z,\lambda)p_+(z,\lambda),
\EEQ
where $\chi\in\Lambda\LieSU(2)_\sigma$ and
$p_+(z)\in\Lambda^+\LieSL(2,\Bcc)_\sigma$. By setting $z=0$
in~\bref{gmhut} and using Eq.~\bref{initialx}, we get
\BEQ
I=\chi(\lambda)p_+(0,\lambda),
\EEQ
which implies that 
\BEQ
\chi(\lambda)=p_+\inv(0,\lambda)\in\LieU(1)
=\Lambda\LieSU(2)_\sigma\cap\Lambda^+\LieSL(2,\Bcc)_\sigma
\EEQ
is $\lambda$-independent. Then $\chi
g_-\chi\inv\in\Lambda^-_\ast\LieSL(2,\Bcc)_\sigma$ and 
the uniqueness of the Birkhoff decomposition gives~\bref{conclusion}.

b)$\Leftrightarrow$c) is trivial.

b)$\Rightarrow$d) Via Iwasawa splitting we obtain
\BEQ
\hF=\chi F\chi\inv
\EEQ
for the extended frames, and thus by Sym's formula:
\BEQ
J(\hPsi_\lambda)=\chi J(\Psi_\lambda)\chi\inv.
\EEQ

d)$\Rightarrow$a) 
Since the surfaces coincide up to a Euclidean motion, their metrics
$\frac{1}{2}e^{u}\diffz\diffzbar$ and $\frac{1}{2}e^{\hu}\diffz\diffzbar$ 
coincide. This implies Eq.~\bref{conditionb} by Eq.~\bref{genDhD}.
\QED

Finally, we consider a case of special interest to us.

\myprop{} {\em 
Consider two CMC-immersions $\Psi:\FD\rightarrow\threespace$ and
$\hPsi:\FD\rightarrow\threespace$ with associated meromorphic
potentials $\xi$ and $\hxi$. Assume moreover, that the metrics induced
by $\Psi$ and $\hPsi$ are complete and that there exists some
proper Euclidean motion $\tT$, s.t.\ $\tT\Psi(\FD)=\hPsi(\FD)$. Then 
\begin{description}
\item[a)] there exists a biholomorphic automorphism $g\in\Aut\FD$,
s.t.\ the diagram

\begin{center}
\setlength{\unitlength}{0.0005in}%
\begin{picture}(2681,1680)(3629,-3007)
\put(6001,-2836){\makebox(0,0)[b]{$\threespace$}}
\put(3676,-2236){\makebox(0,0)[b]{$\Psi$}}
\put(6226,-2236){\makebox(0,0)[b]{$\hPsi$}}
\put(4951,-2611){\makebox(0,0)[b]{$\tT$}}
\put(4951,-1411){\makebox(0,0)[b]{$g$}}
\put(3901,-2836){\makebox(0,0)[b]{$\threespace$}}
\thicklines
\put(4201,-1561){\vector( 1, 0){1500}}
\put(4201,-2761){\vector( 1, 0){1500}}
\put(3901,-1861){\vector( 0,-1){600}}
\put(6001,-1861){\vector( 0,-1){600}}
\put(3901,-1636){\makebox(0,0)[b]{$\FD$}}
\put(6001,-1636){\makebox(0,0)[b]{$\FD$}}
\end{picture}
\end{center} 

commutes,
\item[b)] for the maps $\Psi$ and $\hPsi\circ g$, the associated
meromorphic potentials $\xi$ and $g_\ast\hxi$, and the coefficients
for the Hopf differentials $E=(\hE\circ g)(g^\prime)^2$, the statements of
Theorem~\ref{aundb} are valid.
\end{description}
}

\Proof 
a) Choose any $z_0,\hz_0\in\FD$, satisfying
$\hPsi(\hz_0)=\tT\Psi(z_0)$. Then for sufficiently small neighbourhoods
$U$ and $\hU$ of $z_0$ and $\hz_0$, respectively, $\tT$ induces an
isometry from $U$ to $\hU$. Therefore, from \cite[Sect.~I.11]{He:1},
existence of $g$ follows.

b) clear.
\QED

\newsection \label{summaryth}
It remains to apply Theorem~\ref{aundb} and Corollary~\ref{aundb} to
the case of symmetric CMC-immersions.
First we treat the general case:

\mytheorem{} {\em
Let $\Psi:\FD\longrightarrow\threespace$ be a CMC-immersion with
meromorphic potential 
$$\xi=\lambda\inv\tmatrix0f{E\over f}0\diffz.$$ Let
$g\in\Aut\FD$ be such that $(E\circ g)\cdot(g^\prime)^2=E$.
Let furthermore $g_-$ be the solution of Eq.~\bref{vierzehn2} for $\xi$.
Then the following are equivalent:

1. The immersions $\Psi$ and $\Psi\circ g$ are related by a proper 
Euclidean motion in $\threespace$, i.e., $g\in\Aut_\Psi\FD$. 

2. The surfaces $(\FD,\Psi)$ and
$(\FD,\Psi\circ g)$ have the same conformal metric.

3. There exists a
unitary $z$-independent $\chi(g,\lambda)\in\Lambda\LieSU(2)_\sigma$, s.t.\
$$
g_-(g(z),\lambda)=\chi(g,\lambda)g_-(z,\lambda)p_+(g,z,\lambda),
$$
with a meromorphic function $p_+$ taking values in
$\Lambda^+\LieSL(2,\Bcc)_\sigma$.

4. The equation
\BEQ \label{fourthcondition}
\left|{\hf\over f}\right| = {((\gquer g_-)_{+0})_{11}\over
((\hgquer\hg_-)_{+0})_{11}}
\EEQ
is satisfied, where $\hg_-=g_-\circ g$, $\hf=(f\circ g)g^\prime$
and $(\cdot)_{11}$ denotes the upper left entry of a two by two matrix.
}

\Proof
If we define $\hxi=g_\ast\xi$, $\hE=(E\circ g)(g^\prime)^2$, 
and $\hPsi=\Psi\circ
g$, then $\hg_-$ and $g_-$ satisfy Eq.~\bref{arbitraryinitial} with $g_0=I$,
$\hg_0=g_-(g(0),\lambda)$. 
The equivalence 1.$\Leftrightarrow$3.$\Leftrightarrow$4.\ follows
immediately from Theorem~\ref{aundb}. 
1.$\Rightarrow$2.\ is trivial. 2.$\Rightarrow$1.\ follows
from $\hE=E$ and the second part of Corollary~\ref{DGL}, since both,
$\Psi$ and $\hPsi$, are CMC-surfaces with the same mean curvature.
\QED

Corollary~\ref{aundb} translates into

\mycorollary{} {\em
Let $\Psi$, $\xi$ and $g_-$ be defined as in Theorem~\ref{summaryth}. Let
$g\in\Aut\FD$ be s.t.\ $(E\circ g)(g^\prime)^2=E$. We further assume,
that $g_-(g(0),\lambda)=I$ for all $\lambda\in S^1$.
Then the following statements are equivalent:

a) For $\hf=(f\circ g)g^\prime$ and $\hg_-=g_-\circ g$,
Eq.~\bref{fourthcondition} is satisfied.

b) There exists a $\lambda$-independent, unitary matrix 
$\chi\in\LieU(1)$, s.t.\
\BEQ \label{apart42}
g_-(g(z),\lambda)=\chi g_-(z,\lambda)\chi\inv.
\EEQ

c) There exists a constant $c\in\Bcc$, $|c|=1$, s.t.\ $(f\circ g)g^\prime=f$.

d) There exists a $\lambda$-independent, unitary matrix
$\chi\in\LieU(1)$, s.t.\ for the associated family $\Psi_\lambda$ we
have in the spinor representation
\BEQ
J(\Psi_\lambda\circ g)=\chi J(\Psi_\lambda)\chi\inv.
\EEQ

e) $g\in\Aut_\Psi\FD=\Aut_{\Psi_\lambda}\FD$ and
$\psi_\lambda(g)=\psi(g)$ for all $\lambda\in S^1$.

f) The immersions $\Psi_\lambda$ and $\Psi_\lambda\circ g$ are for all
$\lambda\in S^1$ related by the same rotation around the $e_3$-axis in 
$\threespace$.
}

\Proof
As in the proof of Theorem~\ref{summaryth} we set $\hg_-=g_-\circ g$,
$\hxi=\hg_-\inv\diff\hg_-=g_\ast\xi$, $\hPsi=\Psi\circ g$, and
$\hE=(E\circ g)(g^\prime)^2$.  By assumption, $\hE=E$ and
$\hg_-(0,\lambda)=g_-(0,\lambda)=I$. Therefore, with these
definitions, the assumptions of Corollary~\ref{aundb} are satisfied,
and the equivalence a)$\Leftrightarrow$b)$\Leftrightarrow$c)$\Leftrightarrow$d)
follows immediately from Corollary~\ref{aundb}. 

d)$\Rightarrow$f)
{}From d) it follows (see also the remark
after Eq.~\bref{Symtrafo}), that $\Psi_\lambda$ and
$\Psi_\lambda\circ g$ are related by a $\lambda$-independent rotation $R$
around the $e_3$-axis in $\threespace$. 

f)$\Rightarrow$e) By f),
$g\in\Aut_{\Psi_\lambda}\FD$ for all $\lambda$. If we define
$\psi_\lambda:\Aut_{\Psi_\lambda}\FD\rightarrow\Aut\Psi_\lambda(\FD)$
as in Section~\ref{Symsection}, then $R=\psi_\lambda(g)=\psi(g)$ 
for all $\lambda\in S^1$.

e)$\Rightarrow$d) By Theorem~\ref{loopinterpretation}, to
$g\in\Aut_\Psi\FD$ there is associated a matrix 
$\chi(g,\lambda)\in\Lambda\LieSU(2)_\sigma$. 
If we set $\tT_\lambda=(R_{\tT_\lambda},t_{\tT_\lambda})=\psi_\lambda(g)$,
then we see by Eq.~\bref{rottrans},
that for every $\lambda$, $\tT_\lambda$ determines $\chi(g,\lambda)$
up to a sign. Since $\chi(g,\lambda)$ depends continuously on
$\lambda$, we have, that with $\tT_\lambda$ also $\chi(g,\lambda)$ is
independent of $\lambda$. Therefore, $\chi(g,\lambda)$ is a
diagonal, unitary matrix, from which, by Eq.~\bref{Symtrafo}, d) follows.
\QED

\newsection
In the next two sections we want to summarize the previous
considerations and, as a conclusion, we want to provide a recipe 
for the construction of
CMC-immersions of a fixed Riemann surface $M$.

Let us start with $\FD$ the complex plane or the unit disk, and let us
choose a group $\Gamma$ of biholomorpic automorphisms of $\FD$, which
acts freely and discontinuously. Then $M=\Gamma\setminus\FD$ is a
Riemann surface~\cite{FaKr:1,Spr:1}.

Let us further assume, that there exists a CMC-immersion
$\Phi:M\rightarrow\threespace$, s.t.\ $M$ with the induced metric is
complete. We obtain the following commutative diagram of Section~\ref{DGL}:

\begin{center}
\setlength{\unitlength}{0.00083300in}%
\begin{picture}(1265,1215)(4316,-2269)
\thinlines
\put(4801,-1261){\vector(-1,-2){300}}
\put(4951,-1261){\vector( 1,-2){300}}
\put(4501,-2086){\makebox(0,0)[b]{$M$}}
\put(4651,-2011){\vector( 1, 0){450}}
\put(5251,-2086){\makebox(0,0)[b]{$\threespace$}}
\put(5251,-1561){\makebox(0,0)[b]{$\Psi$}}
\put(4801,-2236){\makebox(0,0)[b]{$\Phi$}}
\put(4876,-1186){\makebox(0,0)[b]{$\FD$}}
\put(4501,-1561){\makebox(0,0)[b]{$\pi$}}
\end{picture}
\end{center}

We can assume, that $\Psi$ is a conformal immersion. 

We consider the frame $\hF:\FD\rightarrow\LieSO(3)$ given by 
\BEQ
\hF=(e^{-\frac{u}{2}}\Psi_x,e^{-\frac{u}{2}}\Psi_y,
\frac{\Psi_x\times\Psi_y}{|\Psi_x\times\Psi_y|}).
\EEQ
We normalize $\Phi$ and $\Psi$ by assuming $\hF(0)=I$. The frame $\hF$
can be lifted to a frame $\tF:\FD\rightarrow\LieSU(2)$.

Considering the differential $\talpha=\tF\inv\diff\tF$, we decompose
$\talpha=\talpha_p+\talpha_k$, where $\talpha_k$ is the diagonal, and
$\talpha_p$ is the off-diagonal part of $\talpha$. Then $\talpha_p$
decomposes into a $\diffz$- and a $\diffzbar$-part:
$\talpha_p=\talpha_p^\prime\diffz+\talpha_p^\dprime\diffzbar$.

We set
\BEQ \label{alphaform5}
\alpha^\lambda=\lambda\inv\talpha_p^\prime\diffz+\talpha_k+\talpha_p^\dprime
\diffzbar.
\EEQ
We know \cite{DoPeWu:1}, that the Gauss map of $\Phi$ is
harmonic iff $\alpha^\lambda$ is integrable, i.e., if there exists an
extended frame 
\BEQ
F:\FD\times S^1\rightarrow\LieSU(2),\kern1cm F=F(z,\zquer,\lambda),
\EEQ
s.t.\ 
\BEQ \label{Fdef5}
F\inv\diff F=\alpha^\lambda,\kern2cm F(0,0,\lambda)=I.
\EEQ
Then, of course, $F(z,\zquer,\lambda=1)=\tF(z,\zquer)$.

We note, that for every fixed $\lambda$, Sym's
formula~\cite{Sym:1,Bo:1} gives a CMC-immersion $\Psi_\lambda$:
\BEQ
\Psi_\lambda(z,\zquer)=\partial_\theta F\cdot F\inv+\frac{i}{2}F\sigma_3
F\inv,\kern0.5cm\lambda=e^{i\theta}.
\EEQ
For $g\in\Gamma$, we have (Theorem~\ref{loopinterpretation})
\BEQ \label{Fg5}
(F\circ
g)(z,\zquer,\lambda)=\chi(g,\lambda)F(z,\zquer,\lambda)k(g,z,\zquer).
\EEQ
where $\chi:\Gamma\rightarrow\Lambda\LieSU(2)_\sigma$ is a homomorphism
up to a sign, and $k$ is a cocycle up to a sign.
Projecting $F$ back into $\LieSO(3)$ then actually gives a
homomorphism and a cocycle, respectively.

{}From~\bref{Fg5} we know, that $\chi(g,\lambda=1)=\pm I$ for
$g\in\Gamma$.

Conversely, let $F$ be an extended frame, i.e., $F$ satisfies
Eq.~\bref{Fdef5} with $\alpha^\lambda$ of the form~\bref{alphaform5}.
Assume now, that $F$ satisfies Eq.~\bref{Fg5} for all
$g\in\Gamma$. Then, for every fixed $\lambda$, Sym's formula gives a
CMC-immersion $(\FD,\Psi_\lambda)$, s.t.\
$\Gamma\subset\Aut_{\Psi_\lambda}\FD$.

If there exists some $\lambda_0\in S^1$, s.t.\ $\chi(g,\lambda_0)=\pm
I$ for each $g\in\Gamma$, i.e., such that $\Psi_{\lambda_0}\circ
g=\Psi_{\lambda_0}$ for all $g\in\Gamma$, then $\Psi_{\lambda_0}$
factors through $M=\Gamma\setminus\FD$, i.e., there exists a
CMC-immersion $\Phi:M\rightarrow\threespace$.

{\em Therefore, in order to
find CMC-immersions of $M=\Gamma\setminus\FD$, it suffices to find an
extended frame satisfying Eq~\bref{Fg5}, s.t.\ for some
$\lambda_0\in S^1$, the kernel equation $\chi(g,\lambda_0)=\pm I$ holds for
all $g\in\Gamma$.}

\newsection
Assume now that the extended frame $F$ satisfies
Eq.~\bref{Fg5}.
We split $F=g_-g_+$ for $z\in\FD\setminus\FS$, where $\FS\subset\FD$
is the discrete set of points where the Birkhoff splitting of
$F$ fails. For 
$\xi=g_-\inv\diff g_-$ we then obtain:
\BEA
\xi & = & \lambda\inv\tmatrix0f{E\over f}0\diffz, \label{cond1}\\
g_\ast\xi & = & p_+\inv\xi p_++p_+\inv\diff p_+,\kern1cm g\in\Gamma,
\label{cond2}
\EEA
where $p_+=p_+(g,z,\lambda)$ has no negative $\lambda$-coefficients
and is holomorphic on $\FD\setminus\FS$.
Moreover, $p_+$ is uniquely defined by $g$ up to a sign.
We note, that the conditions \bref{cond1} and \bref{cond2} imply
for $g_-$, defined by 
\BEQ \label{gmin5}
g_-\inv\diff g_-=\xi, \kern2cm g_-(0,\lambda)=I,
\EEQ
the relation 
\BEQ \label{chi5}
g_-\circ g=\chi g_-p_+,
\EEQ
where, for each $g\in\Gamma$, $\chi=\chi(g,\lambda)$ is
$z$-independent, see Section~\ref{automorphism}.
{}From this we cannot conclude the identity~\bref{Fg5}. This is the
case, if actually $\chi(g,\lambda)\in\LieSU(2)$ for all
$g\in\Gamma$, $\lambda\in S^1$. Below we give conditions that imply 
$\chi(g,\lambda)\in\LieSU(2)$.

Evaluating the $\lambda\inv$-coefficient of \bref{cond2} we obtain
(see Section~\ref{automorphism})
\BEA \label{Efcond1}
(E\circ g)(z)(g^\prime(z))^2 & = & E(z),\\
(f\circ g)(z)g^\prime(z) & = & r^{-2}(g,z)f.\label{Efcond2}
\EEA
Here $r^{-2}$ satisfies the cocycle condition
\BEQ
r^{-2}(g_2\circ g_1,z)=r^{-2}(g_2,g_1(z))r^{-2}(g_1,z),\;g_1,g_2\in\Gamma.
\EEQ
Moreover, if $\chi(g,\lambda)$ in~\bref{chi5} is unitary for all
$\lambda\in S^1$, then (Theorem~\ref{summaryth}) 
$|r(g,z)|$, $g\in\Gamma$, is explicitly given by
\BEQ\label{Efcond4}
|r(g,z)|^2= {((\overline{g_-\circ g}^\top (g_-\circ g))_{+0})_{11}\over
((\gquer g_-)_{+0})_{11}}.
\EEQ
Conversely, let $E$ be holomorphic on $\FD$ and let $f$ be meromorphic on
$\FD$, s.t.\ 
$$
\xi=\lambda\inv\tmatrix0f{E\over f}0\diffz
$$
defines a smooth CMC-immersion $(\FD,\Psi)$ with associated family
$\Psi_\lambda$ (see~\cite{DoHa:1}). {\em If $E$ and $f$ satisfy 
Eqs.~\bref{Efcond1},~\bref{Efcond2}, and~\bref{Efcond4} for all
$g\in\Gamma\subset\Aut\FD$, then (Theorem~\ref{summaryth})
$\Gamma\subset\Aut_\Psi\FD=\Aut_{\Psi_\lambda}\FD$.}
In particular, the solution $g_-$ of~\bref{gmin5} satisfies Eq.~\bref{chi5}
with $\chi=\chi(g,\lambda)\in\Lambda\LieSU(2)_\sigma$ and
$p_+\in\Lambda^+\LieSL(2,\Bcc)_\sigma$. 
{\em If we want to factor $\Psi_{\lambda_0}$ for some $\lambda_0\in S^1$ 
through $M=\Gamma\setminus\FD$,
then we need, in addition, a fourth condition
\BEQ \label{Efcond5}
\chi(g,\lambda_0)=\pm I, \kern1cm\mbox{\rm for all $g\in\Gamma$}.
\EEQ
}

In a brief summary: Every CMC-immersion $(M,\Phi)$ (with the exception
of a sphere), in particular every compact CMC-immersion,
can be obtained by the following construction:

Let $\FD$ be the complex plane or the unit disk. Let $\Gamma$ be a
Fuchsian or elementary group acting on $\FD$. Let $E$
be holomorphic on $\FD$ and $f$ meromorphic on $\FD$, s.t.\ 
$$
\xi=\lambda\inv\tmatrix0f{E\over f}0\diffz
$$
defines a smooth CMC-immersion $(\FD,\Psi)$.

If $E$ and $f$ satisfy Eq.~\bref{Efcond1},~\bref{Efcond2},
and~\bref{Efcond4} for all $g\in\Gamma$, then Eq.~\bref{chi5} holds,
with a unitary matrix $\chi=\chi(g,\lambda)$, for all $g\in\Gamma$.
If, in addition, 
there exists some $\lambda_0\in S^1$, s.t.\ Eq.~\bref{Efcond5}
is satisfied for all $g\in\Gamma$, then $(\FD,\Psi_{\lambda_0})$
covers a CMC-immersion $(M,\Phi)$ with Fuchsian group $\Gamma$.

It would be interesting to see the program above being carried out for
compact CMC-surfaces explicitly.

\section{Two examples} \label{examples}\message{[examples]}

\newsection \label{Smythexample}
As an example for the case treated in Corollary~\ref{summaryth},
we revisit the Smyth
surfaces, where the metric is rotationally symmetric.

First we characterize Smyth surfaces in terms of their meromorphic
potentials.

\myprop{} {\em The Smyth surfaces are, up to coordinate
transformation, in $1-1$-correspondence to the
meromorphic potentials of the form 
\BEQ \label{xiSmythform}
\xi=\lambda\inv\tmatrix01{cz^m}0\diffz,
\EEQ
where $m=0,1,2,\dots$ is an integer and $c\neq0$ a complex number.
}

\Proof Let $(\FD,\Psi)$ be a CMC-immersion that admits a one-parameter
group of rotations as self-isometries. W.l.o.g.\ we assume, that the
center of these rotations is $z_0=0$. From Corollary~\ref{Cautom}
we thus obtain, that up to a reparametrization $E=d z^m$, 
$d\in\Bcc\setminus0$. Moreover, for the function $f$ we have
by Wu's formula~\cite{Wu:1}: $f(z)=\frac{1}{4}e^{u(z,0)-\frac{1}{2}u(0,0)}$. 
Since $u=u(r)$, we have
$u(z,0)=u(0,0)$ and $f=f_0=\const$ follows.
If we change coordinates by $\diff w=f_0\diffz$, then in
$w$-coordinates the meromorphic potential is of the
form~\bref{xiSmythform}.
That all $c$ and all $m$ actually occur, follows from~\cite{Sm:2}.
\QED

\newsection
Consider now $\xi$ as in Theorem~\ref{Smythexample} and the
corresponding CMC-immersion
$\Psi=\Psi_c^m:\Bcc\rightarrow\threespace$,
which has a single umbilic of order $m$ at $z=0$.

We will first describe the associated family of
$(\Bcc,\Psi_\lambda)=(\Bcc,(\Psi_c^m)_\lambda)$ geometrically. 

In the DPW construction we get $\Psi_\mu$ for fixed $\mu\in S^1$ by
substituting $\lambda\rightarrow\lambda\mu$ in the formulas for $g_-$,
$F$ and $\Psi$. If we do this, we get for the meromorphic potential
\BEQ
\xi(z,\lambda)\rightarrow\txi(z,\lambda)=\xi(z,\mu\lambda)
=\lambda\inv\tmatrix0{\mu\inv}{\mu\inv cz^m}0\diffz.
\EEQ
Let us set $\mu=e^{i\theta}$, $\theta\in[0,2\pi)$. 
If we rotate the coordinate system on $\FD$ by $\alpha\theta$,
$\alpha\in\Brr$, i.e., $\tz=e^{-i\alpha\theta}z$, then 
\BEQ
\txi(\tz,\lambda)=\lambda\inv\tmatrix0{e^{i(\alpha-1)\theta}}{
e^{i(\alpha(m+1)-1)\theta}c\tz^m}0\diff\tz.
\EEQ
Therefore, if we choose $\alpha=\frac{2}{m+2}$, then
\BEQ \label{Smythmeropassoc}
\txi(\tz,\lambda)=\chi_\mu\lambda\inv\tmatrix01{c\tz^m}0\diff\tz\chi_\mu\inv,
\EEQ
where conjugation with
\BEQ \label{rotangle}
\chi_\mu=\tmatrix{e^{-i\frac{m}{2m+4}\theta}}00{e^{i\frac{m}{2m+4}\theta}}
\EEQ
describes a rotation in $\threespace$ about an angle $-\frac{m}{m+2}\theta$ 
around the $e_3$ axis.

This gives the following

\mylemma{} {\em
For the CMC-immersions $(\Bcc,\Psi_c^m)$, defined by $E=cz^m$, $f=1$, 
all members of the associated family coincide up
to a rigid rotation and a coordinate transformation.
}

\Proof By Eq.~\bref{Smythmeropassoc} 
the variation of the spectral parameter amounts for the meromorphic potential
to a coordinate transformation on $\FD$ plus a rigid rotation in 
$\threespace$.
Since a rigid rotation also leaves the initial condition
for $g_-$ invariant, we get with the same notation as in
Eq.~\bref{Smythmeropassoc} for fixed $\mu\in S^1$
\BEQ
g_-(z,\mu\lambda)=\chi_\mu g_-(e^{-i\alpha\theta}z,\lambda)\chi_\mu\inv, 
\kern0.5cm\mu=e^{i\theta},
\EEQ
where $\alpha=\frac{2}{m+2}$.
This shows, that $\Psi^m_c(\FD)$ and $(\Psi^m_c)_\mu(\FD)$ are related by a
rigid rotation around the $e_3$-axis in $\threespace$. The rotation
angle in $\threespace$ is given by Eq.~\bref{rotangle} as 
$-\frac{m}{m+2}\theta$.
\QED

\newsection
In this section we determine the group $\Aut_{\Psi_c^m}\FD$.

\myprop{} {\em
If $(\Bcc,\Psi_c^m)$ is a nondegenerate Smyth surface, then:
\BEQ
\Aut\Psi_c^m(\FD)\cong\Aut_{\Psi_c^m}\FD=\{g|g(z)=e^{2\pi i{k\over m+2}}z, k=0,\ldots,m+1\}.
\EEQ
}

\Proof From Theorem~\ref{Smythexample} we know that the Hopf
differential is of the form $E=cz^m\diffz^2$ and the metric factor
$e^u$ only depends on the radius. Therefore, Lemma~\ref{secfund} shows, that
a rotation is in $\Aut_{\Psi_c^m}\Bcc$, iff it leaves $cz^m\diffz^2$
invariant.
It is easy to see, that these are exactly the rotations $z\mapsto
e^{\frac{2\pi i}{m+2}}z$. Using Lemma~\ref{SmythDel2}, we see that for
nondegenerate Smyth surfaces we have
$\Aut_{\Psi_c^m}\FD\cong\Aut\Psi_c^m(\FD)$. Altogether, this proves
the claim.
\QED

\newsection
Finally, we want to determine, which of the surfaces $(\Bcc,\Psi_c^m)$
are nondegenerate, i.e., are not cylinders in $\threespace$. Since $m$
gives the order of the single umbilic of the surface, and since
cylinders do not have any umbilics, we only need to consider the case
$m=0$.

We already now, that for $m=0$, $c=1$ we
get a cylinder. For arbitrary $c\neq 0$ we can write
\BEQ
\xi(z,\lambda)=\lambda\inv\tmatrix01c0\diffz=\Omega\xi_0(\tz,\lambda)
\Omega\inv,
\EEQ
where $\xi_0$ is the meromorphic potential of the standard cylinder
(see Section~\ref{cyldressed}), $\tz=\sqrt{c}z$, and
\BEQ
\Omega=\tmatrix{c^{-\frac{1}{4}}}00{c^{\frac{1}{4}}}.
\EEQ
This shows, that $\xi$ is obtained from the standard cylinder by a
coordinate transformation and subsequent dressing with
$h_+=\Omega$. Corollary~1 of Section~\ref{cyldressed} then shows, that
$\Psi_c^0(\FD)$ is a cylinder iff $\Omega$ is unitary, i.e., iff $|c|=1$.

We have proved the following

\myprop{} {\em The meromorphic potential
$\xi=\lambda\inv\tmatrix01{cz^m}0\diffz$, $c\in\Bcc\setminus\{0\}$,
$m\in\Bcc$, gives a nondegenerate Smyth
surface iff $m=0$ and $|c|\neq1$.
}

\newsection \label{branched}
Other simple examples are provided by surfaces with branchpoints.
An example is the surface with the meromorphic potential
\BEQ
\xi=\lambda\inv\tmatrix0{z-z_0}10\diffz,
\EEQ
where $f(z)=E(z)=z-z_0$ has a simple zero at $z_0\in\Bcc^\ast$. While this
potential is holomorphic at $z=z_0$, it cannot be produced from a
CMC-immersion, since $f$ is not the square of a meromorphic function
(see~\cite[Theorem~2.3]{DoHa:1}). However, it does give a
nonsingular immersion and globally a surface
$\Psi(\FD)$, $\FD=\Bcc$, in $\threespace$ with one branchpoint.
To find the ``correct'' meromorphic potential of the punctured surface
$\Psi(\Bcc\setminus\{z_0\})$, we need to take into account the elementary
group of the surface. 

We define a coordinate $\tw$ on the universal cover of $\cstar$ by 
$w=z-z_0$ and $w=e^{\tw}=q(\tw)$,
$q:\Bcc\rightarrow\Bcc\setminus\{0\}$. We get $\tw=\ln w$,
$\diff\tw={1\over w}\diff w$, $q_\ast(f(z)\diffz)=e^{2\tw}\diff\tw$,
$q_\ast(E(z)\diffz^2)=e^{3\tw}\diff\tw^2$. On the universal cover we have
\BEQ
\txi=\lambda\inv\tmatrix0{e^{2\tw}}{e^{\tw}}0\diff\tw.
\EEQ
The initial condition changes from $g_-(z=0)=I$ to
$\tg_-(\tw_0)=I$, where some $\tw_0$ with $q(\tw_0)=-z_0\neq0$ is fixed. 
In the coordinate $\tw$ on the universal cover $\Bcc$, 
the elementary group of the surface is the
translation group generated by $\tw\rightarrow\tw+2\pi i$.

\end{document}